\documentclass[11pt]{amsart}
\usepackage{amsfonts}
\usepackage{amsmath,amsthm,amssymb}
\makeindex
\usepackage{enumerate}

\newcommand{\lm}{{\operatorname{lm}}}
\newcommand{\hilbert}{{\mathcal H}}
\newcommand{\age}{{\mathcal A}}
\newcommand{\profile}{\varphi}

\newcommand{\R}{\mathbb{R}}
\newcommand{\K}{\mathbb{K}}
\newcommand{\N}{\mathbb{N}}
\newcommand{\C}{\mathbb{C}}
\newcommand{\Z}{\mathbb{Z}}
\newcommand{\Q}{\mathbb{Q}}

\newtheorem{theorem}{Theorem}
\newtheorem{lemma}{Lemma}
\newtheorem{corollary}{{\bf Corollary}}
\newtheorem{proposition}{\noindent {\bf Proposition}}

\newtheorem{remark}[theorem]{Remark}

\theoremstyle{definition}

\newtheorem{problem}{\noindent {\bf Problem}}

\newtheorem{fact}[theorem]{\noindent{\bf Fact}}

\newtheorem{example}[theorem]{Example}
\newtheorem{examples}[theorem]{Examples}

\def\endproof{\hfill {\kern 6pt\penalty 500
\raise -0pt\hbox{\vrule \vbox to5pt {\hrule width 5pt
\vfill\hrule}\vrule}}}

\def\centerpicture #1 by #2 (#3){\leavevmode
        \vbox to #2{
        \hrule width #1 height 0pt depth 0pt
        \vfill
        \special{pictfile #3}}}

\def\GenShuffle#1#2#3{\mathbin{
      \hbox{\vbox{
        \hbox{\vrule
              \hskip#2
              \vrule height#1 width 0pt
               }%
        \hrule\vskip#3}%
             \vbox{
        \hbox{\vrule
              \hskip#2
              \vrule height#1 width 0pt
               \vrule }%
        \hrule\vskip#3}%
}}}

\def\shuffle{\,{\mathchoice{\GenShuffle{5pt}{3.5pt}{1pt}}%
                           {\GenShuffle{4pt}{3pt}{1pt}}%
                           {\GenShuffle{3pt}{2pt}{0.7pt}}%
                           {\GenShuffle{2pt}{1.5pt}{0.5pt}}}\,}%

\title{The profile of relations.}

\author{Maurice Pouzet}
\address{ Probabilit\'es-Combinatoire-Statistique\\
\small Universit\'e Claude-Bernard Lyon1\\
\small Domaine de Gerland\\
\small b\^at.  recherche\\
\small 50 Avenue Tony Garnier\\
\small 69365 Lyon cedex 07, France\\
\small \rm e-mail: pouzet@univ-lyon1.fr\\
\small Fax 33 4 37 28 74 80}
\date{\today}
\begin{document}

\dedicatory{Dedicated to Roland Fra\"{\i}ss\'e, at the occasion of his $86^{th}$
birthday.}


\begin{abstract} Le {\it profil } d'une structure relationelle $R$ est la fonction $\varphi_R$
qui  compte pour chaque entier  $n$ le nombre de ses sous-structures \`a  $n$  \'el\'ements, les
sous-structures isomorphes \'etant identifi\'ees. 
Dans cet expos\'e, je donne quelques exemples, notamment des exemples venant des groupes, et pr\'esente quelques faits frappants concernant  le comportement des profils. 
 J'indique le r\^ole jou\'e par quelques notions  de la th\'eorie de l'ordre et  de la combinatoire (eg belordre, alg\`ebre ordonn\'ee, th\'eor\`eme de Ramsey) dans l'\'etude du profil. Comme illustration, je montre  que le profil d'une structure relationnelle $R$ dont l'\^age est in\'epuisable et de hauteur au plus $\omega(k+1)$ satisfait l'in\'egalit\'e $\varphi_R(n)\leq{ {n+k}\choose {k}}$ pour tout entier $n$. Les recherches en cours  sugg\`erent de voir le profil d'une structure relationnelle $R$ comme la fonction de Hilbert d'une alg\`ebre gradu\'ee associ\'ee \`a $R$; un exemple est {\it l'alg\`ebre  d'un \^age}, invent\'ee  par P.J. Cameron. Je  pr\'esente la solution d'une conjecture de Cameron sur l'int\'egrit\'e de l'alg\`ebre d'un \^age, ainsi que quelques progr\`es r\'ecents faits avec Y.Boudabbous  et N.Thi\'ery sur la conjecture que la s\'erie g\'en\'eratrice associ\'ee \`a un  profil  est une fraction rationnelle lorsque ce profil
est born\'e  par un polyn\^ome (et la structure a un noyau fini). 

 \end{abstract}

 \maketitle

\section{introduction}
This paper is a survey about  the properties of a combinatorial function,  the {\it profile of a relational structure}. I present a collection of results, some old,  going back to 1971, some new, some published, some unpublished, and -in order to give to the reader a flavor of the techniques- I  detail some proofs.
An overview of results is given in Section 2.  It is organized around two striking properties of the profile, namely: \emph{the profile of an infinite  relational structure $R$ is non decreasing  and,  provided  that the arity of relations constituting $R$ is bounded or the kernel of $R$ is finite,  its growth rate is either polynomial or faster than every polynomial}. As observed by P.J.Cameron, the Hilbert function of the algebra of invariants of a finite group is a profile. This suggests that the generating series associated to a profile  could be a rational fraction whenever  the profile is bounded by a polynomial (and the relational structure has a finite kernel). I present a positive solution obtained jointly with N.Thi\'ery for the class of relational structures admitting a finite decomposition into monomorphic components. As an illustration, I present the characterization  of tournaments with polynomial profile obtained with Y.Boudabbous. A graded algebra, the  \emph{age algebra}, was   associated  by P.J Cameron to a relational structure $R$ in such a way that its Hilbert function is the profile of $R$.  This algebra and its role are presented in Section 3. The solution of a conjecture of Cameron on the integrity of the age algebra is also presented. 
These two sections contain very few detailed proofs.  Some combinatorial tools needed for the study of the profile are presented in Section 4. An outline of the proof that the  growth rate of a profile is either polynomial or faster than every polynomial is given (see Theorem \ref{profilepouzet3}). One of its ingredients is proved, namely the fact  that the  profile of a relational structure  $R$ whose age is inexhaustible and of height at most  $\omega(k+1)$ satisfies the inequality $\varphi_R(n)\leq{ {n+k}\choose {k}}$ for every non-negative integer $n$ (Theorem \ref{theoremproved}).  

The work presented here benefited from discussions and collaborations with several colleagues.  I am pleased to thank them. I  am particularly pleased to mention  Y.Boudabbous, M.Sobrani and   N.Thi\'ery. Without them, the paper would have been different.

This paper is an outgrowth of a paper presented at the conference in honor of Claude Benzaken, in Grenoble, in september 2002. Since then, its content was presented to several audiences, eg  in Caracas, Kingston, St Denis de la Reunion, Calgary, Tampere, Sfax and,  last  but not the least,  in Hammamet, March 2006,  at  the annual meeting of the "Soci\'et\'e Math\'ematique de Tunisie".  I thank the organizers for offering me this  opportunity, the incentive to write this paper and for their warm hospitality.

  \section{Overview}

\subsection{Definitions and Simple Examples}\label{defexamples}

A \emph{relational structure} is a realization of a language whose
non-logical symbols are predicates. This is a pair $R:= (E,
(\rho_i)_{i \in I})$ made of a set $E$ and of a family of $m_i$-ary
relations $\rho_i$ on $E$. The set $E$ is the \emph{domain} or \emph{base} of $R$. The family $\mu:=(m_i)_{i\in I}$ is the
\emph{signature} of $R$. The \emph{substructure induced by $R$ on a
subset $A$ of $E$}, simply called the \emph{restriction of $R$ to $A$},
is the relational structure $R_{\restriction A}:= (A, (A^{m_i}\cap
\rho_i)_{i\in I})$. Notions of {\it isomorphism} and {\it local
isomorphism} from a relational structure to an other one are defined in
a natural way as well as the notion of {\it isomorphic type} (see
Section 4 for undefined notions). In the sequel, $\tau(R)$ stands for
the isomorphic type of a relational structure $R$ and $\Omega_{\mu}$
stands for the set of isomorphic types of finite relational structures
with signature $\mu$.

The \emph{profile} of $R$ is the function
$\profile_R$ which counts for every integer $n$ the number
$\profile_R(n)$ of substructures of $R$ induced on the $n$-element
subsets, isomorphic substructures being identified.

Clearly, this
function only depends upon the set $\age(R)$ of finite substructures
of $R$ considered up to an isomorphism, a set introduced by R.~Fra{\"\i}ss{\'e}
under the name of \emph{age} of $R$ (see~\cite{fraissetr}).

If the signature $\mu$ is finite (in the sense that $I$ is finite),
there are only finitely many relational structures with signature $\mu$
on an $n$-element domain, hence $\profile_R(n)$ is necessarily an
integer for each integer $n$. In order to capture examples coming from
algebra and group theory, we
cannot preclude $I$ to be infinite. But then, $\varphi_R(n)$ could be an infinite cardinal.

As far as we will be   concerned by the behavior of $\varphi_R$,  we will exclude this case. Indeed, we have:
\begin{fact}\label{infiniteprofile2}
Let $n<\vert E\vert $. Then
\begin{equation}\label{infiniteprofile}
 \varphi_R(n)\leq (n+1)\varphi _R(n+1)
\end{equation}
In particular:
\begin{equation}
\text{If} \; \varphi_R(n) \; \text{ is infinite then}\;
\varphi_R(n+1) \; \text{ is infinite too and }\; \varphi_R(n)\leq
\varphi _R(n+1).
\end{equation}
\end{fact}

Inequality (\ref{infiniteprofile}) follows from a simple counting
argument. Let $[E]^n$ be the set of $n$-element subsets of $E$,
$\age(R)_n: =\{\tau(R_{\restriction F}): F\in [E]^n\}$ and $[E]^{n+1}$,
$\age (R)_{n+1}$ be the sets defined in a similar way. Let $\Gamma:=
\{(\tau (R_{\restriction F}), \tau(R_{\restriction F\cup \{x\}})): F\in
[E]^n\; \text{and}\; x\in E\setminus F\}$. We have trivially
$\varphi_R(n)\leq \sum_{\tau\in \age(R)_n}\vert \{(\tau, \nu)\in \Gamma
\}\vert=\vert \Gamma\vert= \sum_{\nu\in \age(R)_{n+1}}\vert \{(\tau,
\nu)\in \Gamma\} \vert \leq (n+1) \varphi_R(n+1)$. Item (2) follows
immediately.

Hence, except in very few occasions, mentioned in the text,
we  make the assumption that
$\profile_R$ is integer valued, no matter how large $I$ is.
With this assumption, \emph{profiles of relational structures with
bounded signature are profiles of relational structures with finite
signature}, structures that R. Fra{\"\i}ss{\'e} call
\emph{multirelations}, a fact that we record for further use.

\begin{fact}\label{finiteprofile} 
Let $R:= (E,
(\rho_i)_{i \in I})$ be a relational structure of signature $\mu:=
(m_i)_{i\in I}$ and let $n$ be a non-negative integer.
\begin{enumerate}
\item If $\varphi_R(n) $ is an integer there is a \emph{finite} subset
$I'$ of $I$ (whose size is bounded by a function of $n$ and
$\varphi_R(n)$) such that the local isomorphisms of $R$ and of its
reduct $R^{\restriction I'}:= (E, (\rho_{i})_{i\in I'})$ which are
defined on the $n$-element subsets of $E$ are the same, hence
$\varphi_R(n)= \varphi_{R^{\restriction I'}}(n)$.
\item If in addition $\mu$ is bounded above by $n$, that is $\max{\mu}:=
\max \{m_i: i\in I\} \leq n$, one may choose $I'$ such that $\varphi_R=
\varphi_{R^{\restriction I'}}$.
\end{enumerate}
\end{fact}

Several counting functions are profiles.
Here is some simple minded examples.
\begin{enumerate}
\item {\it The binomial  coefficient } ${n+k \choose k}$.
 Let  $R:=(\Q, \leq, u_{1},\dots,u_{k})$ where $\leq $ is the natural order on the set $\Q$  of rational numbers,
  $u_{1},\dots,u_{k}$ are $k$ unary relations which divide $\Q$ into $k+1$ intervals. Then
$\varphi_{R}(n) = {n+k \choose k}$.

\item  {\it The exponential} $n \hookrightarrow k^{n}$.  Let $R:=(\Q, \leq , u_1, \dots, u_k)$, where again
  $u_{1}, \dots, u_{k}$ are $k$ unary relations, but  which divide $\Q$ into $k$ ``colors''
  in such a way that between two rational numbers all colors appear. Then $\varphi_{R}(n) = k^n$.

\item {\it The factorial} $n \hookrightarrow n!$. Let $R:= (\Q, \leq,
\leq')$, where $\leq'$ is an other linear order on $\Q$ such a way that
the finite restrictions induce all possible pairs of two linear orders
on a finite set (eg take for $\leq'$ an order with the same type as the
natural order on the set $\N$ of non-negative integers). Then
$\varphi_{R}(n) = n!$

\item  {\it The partition function} which counts  the number $p(n)$ of partitions of the integer
$n$.   Let $R:= (\N, \rho)$ be the infinite  path on the integers  whose edges  are pairs $\{x,y\}$ such that $y = x+1$. Then $\varphi_R(n)=p(n)$. The determination
of its asymptotic growth is a famous achievement, the difficulties encountered to prove that
$p(n) \simeq
\frac{1}{4
\sqrt{3n}} e^{\pi \sqrt{\frac{2n}{3}}}$ (Hardy and Ramanujan,  1918) suggest some difficulties in the
general study of profiles.
\end{enumerate}

An important class of functions comes from permutation groups.  The {\it orbital profile} of a permutation
group  $G$ acting on a set $E$ is the function $\theta_{G}$ which counts for each integer
$n$
the number, possibly infinite, of orbits of the  $n$-element subsets of $E$.  As it is easy to see, $\theta_{G}$ is the profile of some relational structure $R:= (E,
(\rho_i)_{i \in I})$ on $E$. In fact, as it is easy to see:

\begin{lemma}\label{adherence}
For every permutation group  $G$ acting on  a set  $E$ there is a relational structure $R$ on
$E$ such that:
\begin{enumerate}
\item  Every isomorphism $f$  from a finite restriction of  $R$ onto an other
extends to an automorphism of $R$.
\item  $Aut R = \overline {G}$ where $\overline {G}$ is the topological adherence of
$G$ into the symmetric group  $\mathfrak{G}(E)$, equipped with the topology induced by the
product topology on $E^E$, $E$ being equipped with the discrete topology.
\end{enumerate}
\end{lemma}

Structures satisfying condition 1) are called \emph {homogeneous} (or
\emph{ultrahomogeneous}). They are now considered as one of the basic
objects of model theory. Ages of such structures are called {\it
Fra{\"\i}ss{\'e} classes} after their characterization by
R.Fra{\"\i}ss{\'e}\cite{fraisse54}. In many cases, $I$ is infinite, even
if $\theta_G(n)$ is finite. Groups for which $\theta_G(n)$ is always
finite are said {\it oligomorphic} by P.J.Cameron. The study of their
profile is whole subject by itself \cite {cameronbook}. Their relevance
to model theory stems from the following result of Ryll-Nardzewski
\cite{ryll}.

\begin{theorem}
\label{thm}
Let  $G$ acting on a denumerable set $E$ and $R$ be a relational structure such that $Aut R=
\overline{G}$. Then
$G$ is oligomorphic if and only if the complete theory of $R$ is  $\aleph_{0}$-categorical.
\end{theorem}

\subsection{A Sample of Results}

\subsubsection{\bf The Profile Grows}

Inequality (1) given in the previous subsection can be substantially improved:

\begin{theorem}
\label{increaseinfinite}
  If $R$ is a relational structure on an infinite set then
  $\profile_R$ is non-decreasing.
\end{theorem}
This result was conjectured with R.Fra{\"\i}ss{\'e} \cite{fraissepouzet1}. We proved it in $1971$; the proof - for a single relation- appeared in 1971 in R.Fra{\"\i}ss{\'e}'s book ~\cite{fraisseclmt1}, Exercise~8
p.~113;  the general case was detailed in  \cite{pouzetrpe}. The proof relies on Ramsey theorem \cite{ramsey}. We give it in \ref{proofincreaseprofile}.

More is true:
\begin{theorem}\label{linalg}
If $R$ is a relational structure on a set $E$ having  at least  $2n+k$ elements
then
$\varphi_{R}(n) \leq \varphi_{R} (n+k)$.
\end{theorem}

Meaning that if $\vert E\vert :=m$ then  $\varphi_{R}$ increases up to  $\frac{m}{2}$; and,
for
$n
\geq
\frac{m}{2}$ the value in  $n$ is at least the value of the symmetric of  $n$ w.r.t.
$\frac{m}{2}$.\\

The result is a straightforward consequence of   the following property of incidence matrices.

Let $n, k, m$ be three non-negative integers and $E$ be an $m$-element
set. Let $M_{n, n+k}$ be the matrix whose rows are indexed by the
$n$-element subsets $P$ of $E$ and columns by the $n+k$-element subsets
$Q$ of $E$, the coefficient $a_{P,Q}$ being equal to $1$ if $P\subseteq
Q$ and equal to $0$ otherwise.

\begin{theorem}\label {kantor} 
If $2n+k\leq m$ then $M_{n, n+k}$ has full row rank (over the the field of rational numbers).
\end{theorem}

With this result the proof of Theorem \ref{linalg} goes as follows:

We suppose that $\varphi_{R} (n+k)$ is finite (otherwise, from Fact
\ref{infiniteprofile2}, the stated inequality holds). Thus, we may
suppose also that $E$ is finite (otherwise, for each isomorphic type
$\tau$ of $n+k$-element restriction of $R$ we select a subset $Q$ of $E$
such that $R_{\restriction Q}$ has type $\tau$ and we replace $E$ by the
union of the $Q's$). We consider the matrix whose rows are indexed by
the isomorphic types $\tau$ of the restrictions of $R$ to the
$n$-element subsets of $E$ and columns by the $n$-element subsets $P$ of
$E$, the coefficient $a_{\tau, P}$ being equal to $1$ if
$R_{\restriction P}$ has type $\tau$ and equal to $0$ otherwise.
Trivially, this matrix has full row rank, hence if we multiply it (from
the left) with $M_{n, n+k}$ the resulting matrix has full row rank.
Thus, there are $\varphi(n)$ linearly independent colums. These columns
being distinct, the restrictions of $R$ to the corresponding
$(n+k)$-element subsets have diff
 erent isomorphic types,  hence  $\varphi_R(n)\leq \varphi_R(n+k)$.

We proved Theorem \ref{linalg} in 1976 \cite{pouzet76}. The same
conclusion was obtained first for orbits of finite permutation groups by
Livingstone and Wagner,
1965 \cite{livingstone}, and extended to arbitrary permutation groups by
Cameron, 1976 \cite{cameron76}. His proof uses the dual version of
Theorem \ref{kantor}. Later on, he discovered a nice translation in
terms of his age algebra. We present it in \ref{behavior}.

Theorem \ref {kantor} is in W.Kantor 1972 \cite{kantor}, with similar results for affine and vector
subspaces of a vector space.   Over the last 30 years, it as been applied and rediscovered many
times; recently, it was pointed out that it appeared in a 1966 paper of
D.H.Gottlieb \cite{gottlieb}. Nowadays, this is one of the fundamental
tools in algebraic combinatorics. A proof, with a clever argument
leading to further developments, was given by Fra\"{\i}ss\'e in the 1986
edition of his book, Theory of relations, see \cite{fraissetr}.

\subsubsection{\bf Jumps in the Growth of the Profile}\label{jump}
Infinite relational structures with profile constant, equal  to $1$,  
were called {\it monomorphic} and characterized by R. Fra\"{\i}ss\'e who
proved that they where {\it chainable}. Later on, those with profile
bounded, called {\it finimorphic}, were characterized as {\it almost
chainable} \cite{fraissepouzet1}. We present these characterizations in
\ref{subchainability}. Beyond bounded profiles, and
provided that the relational structures satisfy some mild
conditions, there are jumps in the behavior of the profiles: eg. no profile grows as $log\; n$ or $nlog\;  n$.

 Let
$\varphi : \N
\rightarrow \N$ and $\psi : \N \rightarrow \N$. Recall that $\varphi = O(\psi)$ and $\psi$
{\it grows as fast as} $\varphi$ if
$\varphi(n) \leq a \psi
(n)$ for some positive real number  $a$ and  $n$ large enough.  We say that
$\varphi$ and $\psi$ have the
{\it same growth} if $\varphi$ grows as fast as $\psi$ and $\psi$
grows as fast as  $\varphi$. The growth of $\varphi$ is
{\it polynomial of degree $k$}  if $\varphi$ has the same growth as  $n
\hookrightarrow
n^{k}$; in other words there are positive real numbers $a$ and $b$ such
that $an^k\leq \varphi\leq bn^k$ for $n$ large enough. Note that the
growth of $\varphi$ is as fast as every polynomial if and only if
$lim_{n\rightarrow +\infty}\frac{\varphi(n)}{n^{k}}=+\infty$ for every non negative integer $k$.

\begin{theorem}\label{profilpouzet1}
  Let $R := (E, (\rho_i)_{i \in I})$ be a relational structure. The
  growth of $\profile_R$ is either polynomial or as fast as every
  polynomial provided that either the signature $\mu : = (n_i)_{i
    \in I}$ is bounded or the kernel $K(R)$ of $R$ is finite.
\end{theorem}

The \emph{kernel} of $R$ is the set $K(R)$ of $x \in E$ such that
$\age(R_{\restriction E \setminus \{x\}}) \neq \age(R)$.  Relational structures with
empty kernel are those for which their age has the {\it disjoint
embedding property}, meaning that two arbitrary members of the age can
be embedded into a third in such a way that their domains are disjoint
\cite{pouzetrm}. In Fra{\"\i}ss{\'e}'s terminology, ages with the {\it
disjoint embedding property} are said \emph{inexhaustible} and
relational structures whose age is inexhaustible are said
\emph{age-inexhaustible}. We will say that relational structures with
finite kernel are \emph{almost age-inexhaustible}\footnote{In order to
agree with the Fra{\"\i}ss{\'e}'s terminology, we disagree with the
terminology of our papers, in which \emph{inexhaustibility}, resp.
\emph{almost inexhaustibility}, is used for relational structures with
empty, resp. finite, kernel, rather than for their ages.}.

At this point, enough to know that the kernel of any relational
structure which encodes an oligomorphic permutation group is finite
(this fact immediate: if $R$ encodes a permutation group $G$ acting on a
set $E$ then $K(R)$ is the set
union  of the orbits of the
$1$-element subsets of $E$ which are finite. Since  the number of these orbits is at most
$\theta_G(1)$, if $G$ is oligomorphic then  $K(R)$ is  finite).

\begin{corollary}
The orbital profile of an oligomorphic group is either polynomial or faster than every polynomial.
\end{corollary}

Groups with orbital profile equal to $1$ were described by P.Cameron in
1976 \cite{cameron76}. From his characterization, Cameron obtained that
the growth rate of an orbital profile is ultimately constant, or it
grows as fast as a linear function with slope $\frac{1}{2}$.

For groups, and graphs, there is a much more precise result than Theorem \ref{profilpouzet1}. It
is due to Macpherson, 1985 \cite{macpherson85}.

\begin{theorem}\label{thmgrap}
The profile of a graph or a permutation groups grows either as a polynomial or as fast as
$f_{\varepsilon}$, where
$f_{\varepsilon}(n) = e^{n^{\frac{1}{2} - \varepsilon}}$, this for every
$\varepsilon
>0$.
\end{theorem}
Note that the $f_{\varepsilon}$ are somewhat similar to the partition
function. Such growth cannot be prevented. Indeed, the partition
function is the orbital profile of the automorphim group of an
equivalence relation having infinitely many classes, all being infinite.
Such a group is imprimitive. In fact, according to Macpherson 1987
\cite{macpherson87}:
\begin{theorem}\label{thmpri}
If $G$ is primitive then either $\theta_G(n)=1$ for all $n\in \N$, or
$\theta_G(n)>c^n$ for all $n\in N$, where $c:=
2^{\frac{1}{5}}-\epsilon$.
\end{theorem}

Some hypotheses on $R$ are needed in Theorem \ref{profilpouzet1}, indeed
\begin{theorem}
\label{slowprofile}
For every non-decreasing  and unbounded map $\profile : \N \rightarrow \N$,
  there is a relational structure $R$ such that $\profile_R$ is
  unbounded and eventually bounded above by
  $\profile$.
  \end{theorem}
More is true.

Let $f:\N \rightarrow\N$ be a non-decreasing map such that $1\leq
f(n)\leq n+1$ for all $n\in \N$. Let $A:= \{n: f(n')<f(n+1)\; \text{for
all }\; n'<n+1\}$. Let $R:= (\N, (\rho_n)_{n\in A})$ in which each
$\rho_n$ is ${n+1}$-ary, with $(x_1,\dots, x_{n+1})\in \rho_n$ if and
only if $\{x_1,\dots, x_{n+1}\}= \{0, \dots, n\}$. Then $\varphi_R= f$.

The reader will notice that if $f$ is unbounded then the signature of
$R$ is unbounded and also the kernel of $R$ is infinite (equal to $\N$).

The hypothesis about the kernel is not ad
hoc. As it turns out, \emph{if the growth of the profile of a relational
structure with a bounded signature is bounded by a polynomial then its
kernel is finite}(cf.  Theorem \ref{profilepouzet2} and Theorem \ref{profilepouzet3}).

An outline of the proof of Theorem \ref{profilpouzet1} is given in the last section of the paper.

Theorems \ref{profilpouzet1} and \ref{slowprofile} were obtained  in  1978 \cite{pouzettr}.
Theorem \ref{slowprofile} and a part of Theorem \ref{profilpouzet1}
appear in \cite{pouzetrpe}, with a detailed proof showing that the
growth of unbounded profiles of relational structures with bounded
signature is at least linear. The notion of kernel is in
\cite{ pouzettr}, see also \cite{pouzetrm},  \cite{pouzetrpe}, and  \cite{pouzetri} Lemme IV-3.1 p.~37.

\subsection{Polynomial Growth}\label{subsectionpoly}

It is natural to ask:
\begin{problem}\label{growthpoly} 
If the profile of a relational structure $R$ with finite kernel has
polynomial growth, is $\varphi_R(n)\simeq cn^{k'}$ for some positive
real $c$ and some non-negative integer ${k'}$?
\end{problem}

The problem was raised by P.J.Cameron for the special case of orbital
profiles \cite{cameronbook}. Up to now, it is unsolved, even in this
special case.

An example,  pointed out by P.J.Cameron, \cite{cameronbook},  suggests that a stronger property holds.

 Let $G'$ be the wreath product $G':=G\wr
\mathfrak {S} _\omega$ of a permutation group $G$
acting on $\{1,\ldots, k\}$  and of
$\mathfrak{ S}_\omega$,
the symmetric group on $\omega$. Looking at
$G'$   as a permutation group acting on
$E':=\{1,\ldots, k\}\times \omega $, then - as observed by Cameron-
 $\theta_{G'}$ is the {\it Hilbert function} $h_{Inv(G)}$ of the subalgebra $Inv(G)$ of $\C[x_1,\ldots, x_k]$ consisting of polynomials in the $k$
indeterminates $x_1,\ldots, x_k$ which are invariant under the
action of $G$. The value of $ {h}_{Inv(G)}(n)$  is,  by definition,    the dimension $dim(Inv_n(G))$ of the subspace of homogeneous polynomials of degree $n$. As
it is well known, the {\it Hilbert series} of
$Inv (G)$,
$${\mathcal H}(Inv(G), x):=\sum_{n=0}^{\infty} {h}_{Inv(G)}(n)x^n$$ is a  rational fraction of the form
\begin{equation}\label{hilbertpoly}\frac{P(x)}{(1-x)\cdots(1-x^{k})}\end{equation}
  with  $P(0)=1$,
$P(1)>0$,  and all coefficients of $P$ being non negative integers.

\begin{problem}\label{growthgroup} Find an
example of a permutation group $G'$ acting on a
set $E$ with no finite orbit, such that the
orbital profile of $G'$ has polynomial growth,
but is not the Hilbert
function of the invariant ring $Inv(G)$ associated
with a permutation group $G$ acting on a  finite
set.
\end{problem}

 Let us associate  to a relational structure
$R$ whose  profile takes only finite values its  {\it generating series}
$${\mathcal H}_{\varphi_{R}}:= \sum_{n=0}^{\infty} \varphi_{R}(n)x^{n}$$

\begin{problem}\label{growthseries}
 If $R$  has  a finite kernel and $\varphi_R$ is bounded above by some polynomial,
 is the series
    $\hilbert_{\profile_R}$ a rational fraction
    of the form
    \begin{displaymath}\label{quasipoly}
    \frac{P(x)}{(1-x)(1-x^2)\cdots(1-x^k)}
      \end{displaymath}
 with  $P\in \Z[x]$?
\end{problem}

Under the hypothesis above we  do not know if   $\hilbert_{\profile_R}$ is a rational fraction.

It is well known that if a generating function is of the form
$\frac{P(x)}{(1-x)(1-x^2)\cdots(1-x^k)}$ then for $n$ large enough,
$a_n$ is a \emph{quasi-polynomial} of degree $k'$, with $k'\leq k-1$,
that is a
  polynomial $a_{k'}(n)n^{k'}+\cdots+ a_0(n)$ whose coefficients
  $a_{k'}(n), \dots, a_0(n)$ are periodic functions.  Hence, a subproblem is:

\begin{problem}\label{growthquasipoly}
If $R$  has a finite kernel and $\varphi_R$ is bounded above by some polynomial,
is $\profile_R(n) $ a
  \emph{quasi-polynomial} for $n$ large enough?
  \end{problem}

\begin{remark}Since the
  profile is non-decreasing, if $\profile_R(n) $ is a
quasi-polynomial for $n$ large enough then  $a_{k'}(n)$ is eventually
  constant. Hence the profile has polynomial growth in the sense that
$\profile_R(n)\sim c n^{k'}$ for some positive real $c$ and $k'\in \N$. Thus, in this case,
Problem \ref{growthpoly} has a positive solution.
  \end{remark}

In the  theory of  languages,  one of the basic results is that the generating series of a
regular language is a rational fraction (see \cite{berstel}). This
result is not far away from our considerations. Indeed, if $A$ is a
finite alphabet, with say $k$ elements, and $A^*$ is the set of words
over $A$, then each word can be viewed as a finite chain coloured by $k$
colors and $A^*$ can be viewed as the age of the relational structure
made of the chain $\Q$ of rational numbers divided into $k$ colors
in such a way that, between two distinct rational numbers, all colors
appear. This structure was Example (2) in Subsection~1.1.

\begin{problem}
Does the members of the age of a relational structure with polynomial
growth can be coded by words forming a regular language?
\end{problem}

\begin{problem} 
Extend the properties of  regular languages to subsets of  $\Omega_{\mu}$.
\end{problem}

 \subsection{Morphology of Relational Structures  with Polynomial Growth}\label{subsectionmultichain}
 We only have a partial   description of  relational structures  with polynomial growth.

 Let us say that a relational structure $R:= (E,
(\rho_i)_{i \in I})$ is  \emph{almost multichai-nable} if there is a finite subset $F$ of $E$ and an enumeration
$(a_{x, y})_{(x,y)\in V\times L}$ of the elements of $E\setminus F$ by a
set $V\times L$, where $V$ is finite and $L$ is a linearly ordered set,
such that for every local isomorphism $f$ of $L$ the map $(1_V, f)$
extended by the identity on $F$ is a local isomorphism of $R$ (the map
$(1_V, f)$ is defined by $(1_V,f)(x, y):= (x, f(y))$).

\pagebreak

Note that if $L$ is infinite, $K(R)$, the kernel of $R$, is a subset of $F$. Thus we have:
\begin{fact}
An almost multichainable relational structure has a finite kernel.
\end{fact}

The notion of almost multichainability was introduced in \cite{
pouzettr} and appeared in \cite{pouzetri, pouzetrm}. It seems to be a
little bit hard to swallow. The special case $\vert V\vert =1$, for
which the linear order $L$ can be defined on $E\setminus F$, is
discussed in more details in~\ref {subchainability}.

The profile of an almost multichainable relational structure is not
necessarily bounded above by a polynomial (see the last two examples
given in Examples \ref{examplegraphs}).

\begin{problem}\label{growthexponential}   
If the profile of an almost multichainable relational structure is not
bounded above by a polynomial, does his profile has exponential growth?
Is the generating series  a rational fraction?
\end{problem}

\begin{theorem}\label{profilepouzet2} 
If the profile of a relational structure $R$ with bounded signature or
finite kernel is bounded above by a polynomial then $R$ is almost
multichainable.
\end{theorem}
For a proof, see Section 4, Theorem \ref{profilepouzet3}.

There are two cases, in fact opposite cases, for which the profile of an
almost multichainable relational structure is bounded above by a
polynomial.
\begin{enumerate}
\item {\bf Case 1.}  $(1_V, f)$ extended by the identity on $F$   is an automorphism of $R$ for every permutation $f$ of $L$.
\item {\bf Case 2.} For every family $(f_x)_{x\in V}$ of local
isomorphisms of $L$, the map $\cup \{f_x: x\in V\}$ extended by the
identity on $F$ is a local isomorphism of $f$ (the map $\cup\{ f_x: x\in
V\}$ associates $(x, f_{x}(y))$ to $(x,y)$).
\end{enumerate}

A relational structure for which there are $F$ and $(a_{x, y})_{(x,y)\in
V\times L}$ such that Case 1 holds is {\it cellular}. This notion was
introduced by Schmerl \cite{schmerl}. We illustrate it below. Relational
structures for which case 2 holds are illustrated in subsection
\ref{subsection25}.

\subsubsection{The Case of Graphs} 
A {\it directed graph} is a pair
$G:= (E, \rho)$ where $\rho$ is a binary relation on E.
Ordered sets and tournaments are special case of directed graphs. We
will use the term {\it graph} if $\rho$ is irreflexive and symmetric. In
this case $\rho$ is identified with the set $\mathcal E$ of pairs
$\{x,y\}$ of members of $E$ such that $x\rho y$, $G$ is identified with
$(E, \mathcal E)$; the members of $E$ and $\mathcal E$ are the {\it
vertices} and {\it edges} of $G$. We denote by $V(G)$, resp. $E(G)$, the
set of vertices, resp. edges, of $G$.

 In terms of profile,  the class of graphs provides interesting examples.
\begin{examples} \label{examplegraphs}
\begin{enumerate}
\item  $\varphi_{G}(n)$ is constant, equal to  $1$, for every $n\leq \vert V(G)\vert$,   if and only if $\varphi_G(2)\leq 1$, that is  $G$ is a clique or an independent set (trivial).

\item    $\varphi_ G$ is bounded if and only if   $G$ is ``almost constant'' in the Fra{\"\i}ss{\'e}'s  terminology, that is
there is a finite subset  $F_{G}$ of vertices such that two pairs $\{x,y\}$ and $\{x',y'\}$ of
vertices having the same intersection on $F_{G}$ are both edges or both non-edges. This fact is an immediate consequence of  Theorem \ref {finimorphic}.

\item  If $G$ is the direct sum of infinitely many edges, or the direct sum   $K_{\omega}\oplus K_{\omega}$ of two infinite
cliques,  then
$\varphi_{G}(n)=\lfloor \frac{n}{2}
\rfloor +1$, whereas ${\mathcal H}_{\varphi_{G}}=\frac {1}{(1-x)(1-x^2)}$.

\item Let   $G$ be the direct sum $K_{(1, \omega)}\oplus \overline K_{\omega}$ of an infinite wheel and
an infinite independent set, or the direct sum $K_{\omega}\oplus \overline K_{\omega}$ of an
infinite clique and an infinite independent set.  Then
$\varphi_{G}(n)=n$. Hence ${\mathcal H}_{\varphi_{G}}=1+\frac {x}{(1-x)^{2}}$, that we may write
$\frac{1-x-x^2}{(1-x)^{2}}$, as well as $\frac{1+x^3}{(1-x)(1-x^{2})}$.

\item Let  $G$ be  the direct sum of infinitely many $k$-element cliques or the direct sum of $k$ infinite cliques. Then
$\varphi_{G}(n)=p_{k}(n) \simeq
\frac{n^{k-1}} {(k-1)!k!}$ and ${\mathcal H}_{\varphi_{G}}=\frac {1}{(1-x)\cdots (1-x^{k})}$.

\item  If  $G$ is either the direct sum of infinitely many infinite cliques -or an infinite path- then
$\varphi_{G}(n)=p(n)$ the partition function.
\item Let $C:= (E, \leq)$ be a chain and $K_{C, \frac{1}{2}}$ be the graph  whose vertex set is $2\times E$ and the edge set is $ \{\{(0,i), (1, j)\}: i<j\; \text {in}\; C\}$. Such a graph is an \emph{half-complete bipartite graph}.
If $C$ is infinite, then  $2^{n-2} \leq \varphi_{K_{C, \frac{1}{2}}}(n)\leq 2^{n-1}$ \cite{macpherson85}, hence its growth is exponential. In fact, one can check that:
${\mathcal {H}_{K_{C, \frac{1}{2}}}}= \frac{1-2x-x^2+3x^3-x^4}{(1-x)(1-2x)(1-2x^2)}= 1+x+2x^2+3x^3+6x^4+10x^5+20x^6+36x^7+72x^8+136x^9+O(x^{10})$.
\item Let $\tilde{K}_{C, \frac{1}{2}}$ be the graph obtained from $K_{C, \frac{1}{2}}$ by adding all possible edges between vertices of the form $(1, i)$, for $i\in E$. Then $\varphi_{\tilde{K}_{C, \frac{1}{2}}}(n)=2^{n-1}$. \end{enumerate}
\end{examples}

\begin{theorem}
The profile of a graph is bounded by a polynomial if and only if this graph is cellular.
\end{theorem}

A straightforward computation shows that the profile of a cellular graph is bounded by a polynomial.
The converse follows directly from Theorem \ref{profilepouzet2} and
Lemma \ref{cellularlemma} below. A self-contained proof will hopefully
appear in a joint work with S.~Thomass\'e and R.~Woodrow.

\begin{lemma}
\label{cellularlemma} 
The growth of the profile of   almost multichainable graph which is not  cellular is at  least exponential
\end{lemma}
Indeed, let $G$ be an almost multichainable graph. The sets $F$, $V$ and
$L$ which appear in the definition of the almost multichainability of
$G$ satisfy the following conditions: $F, V$ are finite, $V(G)= F\cup
V\times L$ and:
\begin{equation}
\{a, (x, i)\}\in E(G)\;  \text {if and only if}\;  \{a, (x, j)\}\in E(G)
\end{equation}
 \text{for all} \; $a \in F, x\in V, j\in L$
\begin{equation}
\{(x,i), (y,j)\}\in E(G)\;  \text {if and only if}\;  \{(x,i'), (y,j')\}\in E(G)
\end{equation}
\text{for all} \; $x,y\in V, i, j , i',j'\in L$ such that $i\rho j$ and
$i'\rho j'$ where $\rho$ is either the equality relation on $L$ or the
strict order $<$ on $L$.

If $G$ is not cellular then there is some permutation $f$ of $L$ such
that $(1_V, f)$ extended by the identity on $F$ is not an automorphism
of $G$. The map $f$ does not preserve the order on $L$, hence, there
are $i_0,j_0\in L$ and $x,y\in V$ such that $\{(x,i_0), (y,j_0)\}\in
E(G)$ and $\{(x,j_0), (y,i_0)\} \not \in E(G)$.

Let $H:=G_{\restriction \{x,y\}\times L}$. This graph is multichainable,
hence it is entirely determined by the edges belonging to $[\{x,
y\}\times \{i_0,j_0\}]^2\setminus \{(x,j_0), (y, j_0)\}$. There are 16
possible graphs. But, if $L$ is infinite, these graphs yield only two
distinct ages, namely the age of $K_{C, \frac{1}{2}}$ and the age of
${\tilde K} _{C, \frac{1}{2}}$, two graphs described in (7) and (8) of
Examples \ref{examplegraphs}. Hence, they yield at most two distinct
profiles. Their growth rates, as computed in Examples
\ref{examplegraphs}, are exponential, hence the growth rate of
$\varphi_G$ is at least exponential as claimed.

We do not know if Problem \ref{growthpoly} has a positive answer for cellular graphs.
Problem \ref{growthseries} has a positive answer for a special class of
relational structures described in the following subsection.

\subsection{Relational Structures  Admitting a Finite Monomorphic Decomposition}\label{subsection25}

A \emph{monomorphic decomposition} of a relational structure $R$ is a
partition $\mathcal P$ of $E$ into blocks such that for every integer
$n$, the induced structures on two $n$-elements subsets $A$ and $A'$
of $E$ are isomorphic whenever the intersections $ A\cap B$ and $
A'\cap B$ over each block $B$ of $\mathcal P$ have the same size.

This notion was introduced  with N.~Thi\'ery  \cite{pouzetthiery}.

If an infinite relational structure $R$ has a
monomorphic decomposition into finitely many blocks, whereof $k$ are
infinite, then the profile is bounded by some polynomial, whose degree
itself is bounded by $k-1$. Indeed, as one may immediately see:
\begin{equation}
\varphi_R(n)\leq \sum_{s\leq r}{r\choose s}{n+k-1-s\choose k-1}\leq 2^r {n+k-1\choose k-1}
\end{equation}
where $r$ is the cardinality of the union of the finite blocks.

One can say more:
\begin{theorem}[\cite{pouzetthiery}]\label{theorem.quasipolynomial}
  Let $R$ be an infinite relational structure $R$ with a monomorphic
  decomposition into finitely many blocks $(E_i, i\in X)$, $k$ of
  which being infinite. Then, the generating series
  $\hilbert_{\profile_R}$ is a rational fraction of the form:
  \begin{displaymath}
    \frac{P(x)}{(1-x)(1-x^2)\cdots(1-x^k)}.
  \end{displaymath}
\end{theorem}

\begin{corollary}[\cite{pouzetthiery}]
Let $R$ a relational structure as above, then $\varphi_{R}$ has a polynomial
growth and in fact
$\varphi_{R}(n)\sim a n^{k'}$ for some positive real $a$, some non-negative integer $k'$.
\end{corollary}

Recently, with N.Thi\'ery, we proved:

\begin{lemma}
\label{growthdegree} 
If $k$ is the least number of infinite blocks that a monomorphic
decomposition of $R$ may have then $\varphi_{R}(n)\sim a n^{k-1}$ .
\end{lemma}

The proof idea of Theorem \ref{theorem.quasipolynomial} is this. To each
subset $A$ of size $n$ of $E$, we associate the monomial
\begin{displaymath}
  x^{d(A)} := \prod_{i\in X} x_i^{d_i(A)},
\end{displaymath}
where $d_i(A)=|A\cap E_i|$ for all $i$ in $X$. Obviously, $A$ is
isomorphic to $B$ whenever $x^{d(A)}=x^{d(B)}$. The \emph{shape} of a
monomial $x^d=\prod x_i^{d_i}$ is the partition obtained by sorting
decreasingly $(d_i, i\in X)$. We define a total order on monomials by
comparing their shape w.r.t. the degree reverse lexicographic order,
and breaking ties by the usual lexicographic order on monomials w.r.t.
some arbitrary fixed order on $X$. If $\tau\in \age(R)$ its \emph{leading monomial} $\lm (\tau)$ is the
maximum of the monomials $x^{d(A)}$ where $R_{\restriction A}$ has isomorphic type $\tau$.
If  $S$ is subset of $X$, set $x_S:=1$ if $S=\emptyset$ and $x_S:= \prod_{i\in S} x_i$ otherwise. We may write every monomial  $m$ ( distinct from $1$) as  a product $\prod (x_{S_1})^{m_1}\cdots (x_{S_k})^{m_k}$ for a unique sequence $C:=(\emptyset \subsetneq S_1\subsetneq\dots\subsetneq
  S_r \subseteq X)$ of non-empty subsets of $X$. This sequence is the {\it chain support} of  $m$. To such a sequence $C$,  we associate the set $\lm_C$   of leading monomials with chain support $C$.  We prove that one can
realize $\lm_C$ as the linear basis of some ideal of
  a polynomial ring, so that the generating series of $\lm_C$ is
  realized as an Hilbert series. From this,  one concludes easily that
  $\hilbert_{\profile_R}$ has the same form.

The key property of leading monomials  in this proof is
this lemma:

\begin{lemma}[\cite{pouzetthiery}]\label{lemma.addlayer}
  Let $m$ be a leading monomial, and $S\subset X$ be a layer of $m$.
  Then, either $d_i=|E_i|$ for some $i$ in $S$, or $m x_S$ is again a
  leading monomial.
\end{lemma}

The proof of this result relies on
Proposition~\ref{homogeneouscomponent} below for which we introduce
the following definition. Let $R$ be a relational structure on $E$; a
subset $B$ of $E$ is a \emph{monomorphic part} of $R$ if for every
integer $n$ and every pair $A, A' $ of $n$-element subsets of $E$ the
induced structures on $A$ and $A'$ are isomorphic whenever $A\setminus
B=A'\setminus B$.

\begin{proposition}[\cite{pouzetthiery}]\label{homogeneouscomponent}
\begin{enumerate}
\item For every  $x\in E$, the set-union $R(x)$ of all the monomorphic parts
of $R$ containing $x$  is a monomorphic part,
the largest monomorphic part of $R$ containing $x$.
\item The largest monomorphic parts form a monomorphic decomposition of
  $R$ of which every monomorphic decomposition of $R$ is a
  refinement.
  \end{enumerate}
\end{proposition}

We will call {\it canonical} the decomposition of $R$ into maximal
monomorphic parts. This decomposition has the least possible number of
parts.

Despite the apparent simplicity of relational structures admitting a
finite monomorphic decomposition, there are many significant examples.

\subsubsection{Relational Structures which are Categorical for Their Age}

A relational structure $R:= (E, (\rho_{i})_{i
\in I})$ is {\it categorical for its age} if
every $R'$ having the same age than $R$ is
isomorphic to $R$. It was proved in \cite
{hodkinson} that for a relational structure with a
finite signature, this happens just
in case $E$ is at most denumerable  and can be divided into
finitely many blocks such that every permutation
of $E$ which preserves each block is an
automorphism of $R$.

These structures may occur in some interesting areas.  For a simple minded example,  let  $G':=G\wr
\mathfrak{S} _\omega$ be the  wreath product of a  permutation group $G$ acting on $\{1,\ldots, k\}$ and  of
$\mathfrak{S}_\omega$,
the symmetric group on $\omega$. As we have noticed, the orbital profile
of $G'$ is the Hilbert function of the algebra of invariants of $G$.
Now, the group $G'$ is the automorphism group of a relational structure
$R$ which is categorical for its age. Among the
possible $R'$ take
$R':=(E', (\overline \rho_{i})_{i\in I})$ where
$\overline
\rho_{i}:=\{ ((x_{1}, m_{1}),\dots,(x_{n_{i}},
m_{n_{i}})):  (x_{1},\dots, x_{n_{i}})\in
\rho_{i}, ( m_{1}, \dots,
m_{n_{i}})\in \N^{\{1,\dots, n_{i}\}} \}$, $R:=(\{1,\dots,
k\}, (\rho_{i})_{i\in I})$ is  a
relational structure containing the equality
relation and having  signature
$\mu:= (n_{i})_{\in I}$ such that $Aut R= G$.
Such an $R'$ decomposes into $k$ monomorphic
components, namely the sets $E'_{i}, 1\leq i<k$ (where
$E'_{i}:= \{(i,m): m\in \N\}$).

\subsubsection{Quasi-Symmetric Polynomials}

 Let $x_1,\ldots, x_k$ be
$k$-indeterminates and
$n_{1},\dots, n_{l}$ be a sequence of
non-negative integers, $1\leq l\leq k$. The
polynomial
$$\sum _{1\leq i_{1}<\dots<i_{l}\leq
k}x_{i_{1}}^{n_{1}}\dots x_{i_{l}}^{n_{l}}$$ is
 a {\it
quasi-monomial} of degree $n$, where $n=:
{n_{1}}+\dots+{n_{l}}$. The vector space
spanned by the quasi-monomials forms the space
${\mathcal Q}{\mathcal S}_{k}$ of {\it
quasi-polynomials} as introduced by I. Gessel. As
in the example above, the Hilbert series of
${\mathcal Q}{\mathcal S}_{k+1}$ is defined
as
$${\mathcal H}_ {{\mathcal Q}{\mathcal S}_{k}}:=
\sum_{n=0}^{\infty}  dim {\mathcal Q}{\mathcal
S}_{k, n} x^{n}.$$

As shown by F. Bergeron,
C. Retenauer, see \cite{garsia}, this series is  a rational fraction
of the form  $\frac
{P_{k}}{(1-x)(1-x^2)\dots(1-x^{k})}$ where
the coefficients $P_{k}$ are non negative.  Let $R$ be the poset product of a $k$-element
chain by a denumerable antichain. More formally,
$R:= (E, \rho )$ where $E:=
\{1,\ldots, k\}\times \N$ and $\rho:= \{((i,
n), (j, m))\in E$ such that $i\leq j\}$.
Each isomorphic type of an $n$-element
restriction may be identified to a
quasi-polynomial, hence the generating series
associated to the profile of $R$ is
the Hilbert series defined above.
Since $R$
decomposes into $k$ monomorphic
components, the rationality of this series is a
special case of Theorem \ref{theorem.quasipolynomial}.
The reason for which the coefficients of this fraction are non-negative
was elucidated only recently by Garsia and Wallach \cite{garsia}. They
proved that ${\mathcal Q}{\mathcal S}_{k}$ is Cohen-Macaulay.

\subsubsection{Monomorphic Relational Structures, Chainability and Extensions} \label{subchainability}

We present here the origin of the notion of relational structure
admitting a monomorphic decomposition into finitely many blocks.

According to R.Fra{\"\i}ss{\'e} who introduced this notion in 1954 in
his thesis, a relational structure $R:= (E, (\rho_i)_{i \in I})$ for
which $\varphi_R(n)=1$ for every $n\leq \vert E\vert $ is {\it
monomorphic}.

\begin{example}
There are eight kinds of monomorphic directed graphs, four made of
reflexive directed graphs, four made of irreflexive graphs. For the
reflexive ones, there are the chains, the reflexive cliques , the
antichains, plus the $3$-element oriented reflexive cycle. Whereas, for
the irreflexive ones, there are
the acyclic (oriented) graphs, the cliques, the independent sets, and
the $3$-element oriented irreflexive cycle.
\end{example}

Fra{\"\i}ss{\'e} gave a characterization of infinite monomorphic
relational structures by means of his notion of chainability:

A relational structure $R:= (E, (\rho_i)_{i \in I})$ is {\it chainable}
if there is a linear ordering $\leq$ on $E$ such that every local
isomorphism of $L:=(E, \leq)$ is a local isomorphism of $R$.

Since chains are monomorphic, chainable relational structures are also
monomorphic. The converse does not hold, as shown by a $3$-element
cycle. Fra{\"\i}ss{\'e} proved that it holds if the structure is
infinite.

\begin{theorem}\label{chainable}
An infinite relational structure is monomorphic if and only if it is chainable.
\end{theorem}

His proof, given for relational structures of finite signature, was
based on Ramsey's theorem \cite{ramsey} and the compactness theorem of
first order logic. The extension to arbitrary signature requires an
other application of the compactness theorem (for a detailed proof, see
\cite {pouzetrpe}).

We give the proof idea, in the setting of   a generalization  of the monomorphy and chainability notions.

Let $R:= (E, (\rho_i)_{i \in I})$ be a relational structure and $F$ be  a subset of $E$.
The relational structure $R$ is $F$-{\it  monomorphic}  if  for
every non-negative integer $n$   and every  $A,A'\in [E\setminus F]^n$  there is an isomorphism from
$R_{\restriction A}$ onto $R_{\restriction A'}$ which extends by the
identity on $F$ to an isomorphism of $R_{A\cup F}$ onto $R_{A'\cup F'}$.
This relational structure is $F$-{\it chainable} if there is a linear
order $\leq $ on $E\setminus F$ such that every local isomorphism of
$L:= (E\setminus F, \leq)$, once extended by the identity on $F$, is a
local isomorphism of $R$. This relational structure is {\it almost
monomorphic}, resp. {\it almost chainable}, if it is $F$-monomorphic,
resp. $F$-chainable for some finite subset $F$ of $E$.

From Ramsey's theorem, Fra{\"\i}ss{\'e}  deduced the following lemma.

\begin{lemma}
\label{chainablerestriction} 
Let $R$ be a relational structure with domain $E$ and $F$ be a finite
subset of $E$. If the signature of $R$ is finite then there is an
infinite subset $E'$ of $E$ containing $F$ on which
the restriction $R':= R_{\restriction E'}$ is $F$-chainable.
\end{lemma}

Then he applied the compactness theorem of first order logic (in a
weaker form, given by his ``coherence lemma''). Indeed, from Lemma
\ref{chainablerestriction} above, if a monomorphic relational structure
$R$ of finite signature is infinite, it contains an infinite induced
substructure $R'$ which is chainable. Since $R$ is monomorphic, each
finite substructure of $R$ is isomorphic to some finite substructure of
$R'$, hence is chainable. The compactness theorem insures that $R$ is
chainable. As said, this conclusion holds if the signature is arbitrary.

 Theorem \ref{chainable} has the following strenghtening.

 \begin{theorem} \label{monomorphiccomp}
Let $R:= (E, (\rho_i)_{i \in I})$ be a relational structure, $E'$ be a
subset of $E$ and $F:= E'\setminus E$. Let us consider the following
properties
 \begin{enumerate}
 \item[{(i)}]   $R$ is $F$-chainable;
\item[(ii)] $R$ is $F$-monomorphic;
\item[(iii)] $E'$ is a monomorphic part of $R$.
\end{enumerate}
Then $(i)\Rightarrow (ii)\Rightarrow (iii)$. If $E'$ is infinite then
$(ii)\Rightarrow (i)$. If $E'$ is infinite and $E'$ is a monomorphic
component of $R$ then $(iii)\Rightarrow (ii)$.
 \end{theorem}
For these implications, one considers first the case for which the
signature and $F$ are finite. Then, one applies Lemma
\ref{chainablerestriction} and the compactness theorem of first order
logic. The general case follows by another application of the
compactness theorem.

 This yields:
\begin{theorem}\label{finimorphic}
For a relational structure, the following properties are equivalent:
 \begin{enumerate}
\item[{(i)}] The profile of  $R$ is bounded by some integer.
\item[(ii)] $R$ has a monomorphic decomposition into finitely
  many blocks, at most one being infinite.
 \item[(iii)] $R$ is almost chainable.
 \item[(iv)] $R$ is almost monomorphic.
\end{enumerate}

\end{theorem}
Theorem~\ref{finimorphic} above was
proved (without Item (ii))  in~\cite{fraissepouzet1} for finite signature and in~\cite
{pouzetrpe} for arbitrary signature.  Theorem \ref{increaseinfinite} was proved afterward.

From Theorem~\ref{increaseinfinite} and Theorem~\ref{finimorphic}, it
follows that a relational structure $R$ has a monomorphic
decomposition into finitely many blocks, at most one being infinite, if
and only if
\begin{displaymath}
  \hilbert_{\profile_R}=\frac {1+ b_1x+ \dots +b_lx^l} {1-x}\ ,
\end{displaymath}
where $ b_1,\dots, b_l$ are non negative integers.

With this elementary fact in hands, it was not so hard to conjecture the
extension given in Theorem \ref{theorem.quasipolynomial}.

 \subsubsection{An illustration: a proof of Theorem \ref{increaseinfinite}}\label{proofincreaseprofile}

Let $R := (E, (\rho_i)_{i \in I})$ be a relational structure. Suppose
$E$ be infinite. Let $n$ be a non negative integer. We claim that
$\varphi_R (n)\leq \varphi_R(n+1)$.\vskip10pt

\noindent{\bf Case 1.\enskip}$\varphi_R (n)$ is infinite. Then, as stated in
Fact \ref{infiniteprofile2}, $\varphi_R (n)\leq \varphi_R(n+1)$ as
claimed.\vskip6pt

\noindent{\bf Case 2.\enskip}$\varphi_R (n)$ is finite.
We reduce the claim to the case of an almost monomorphic relational structure.\vskip6pt

\noindent{\bf Claim 1.\enskip}There is some finite subset $I'$ of $I$ and some
infinite subset $E'$ of $E$ such that the reduct $R':= R^{\restriction
I'}_{\restriction E'}$ is almost monomorphic and
$\varphi_{R'}(n)=\varphi_R(n)$.\vskip10pt

\noindent{\it Proof of Claim 1.\enskip}According to (1) of Fact \ref
{finiteprofile}, there is some finite subset $I'$ of $I$ such that the
reduct $R^{\restriction I'}:= (E, (\rho_i)_{i \in I'})$ satisfies
$\varphi_{R^{\restriction I'}}(n)=\varphi_R(n)$.
Let $m:= \varphi_{R^{\restriction I'}}(n)$. Select $F_1, \dots F_m$ in
$[E]^n$ such that the restrictions $R^{\restriction I'}_{\restriction
F_1}, \dots, R^{\restriction I'}_{\restriction F_m}$ are pairwise
non-isomorphic. Set $F:= F_1\cup\cdots \cup F_m$. According to Lemma
\ref{chainablerestriction}, there is an infinite subset $E'$ of $E$
containing $F$ such that the restriction $R':= R^{\restriction
I'}_{\restriction E'}$ is $F$-chainable. This restriction is almost
monomorphic. From our construction, $\varphi_{R'}(n)=m$. This proves
Claim~1. \hfill $\blacksquare$\vskip10pt

\noindent{\bf Claim 2.\enskip}If an infinite  relational structure $R':= (E', (\rho_i)_{i\in I'})$ is almost monomorphic then $\varphi_{R'}$ is non-decreasing.

\noindent{\it Proof of Claim 2.\enskip}Let $F$ be a finite subset of $E'$ such that $R'$ is $F$-monomorphic. Let $n$ be a non-negative integer. Let $m:= \varphi_ {R'}(n)$ and let $\tau_1, \dots, \tau_m$ be the isomorphic types of the $n$-element restrictions of $R'$. Select $F_1, \cdots, F_m$ such that  for each $i$, $1\leq i\leq m$, $R'_{\restriction F_i }$ has type $\tau_i$ and $\vert F\cap F_i\vert$ is minimum. Pick $x\in E'\setminus (F\cup F_1\cup \cdots \cup F_m)$ and set $F'_i:=F_i\cup \{x\}$ for $i$, $1\leq i\leq m$.
We claim that the restrictions $R'_{\restriction F'_1 }, \dots,
R'_{\restriction F'_m }$ are pairwise non-isomorphic, from
which the inequality $\varphi_{R'}(n)\leq \varphi_{R'}(n+1)$ will
follow. Indeed, suppose that there is some isomorphism $f$ from
$R'_{\restriction F'_i}$ onto $R'_{\restriction F'_j}$.
With no loss of generality, we may suppose $\vert F_i\cap F\vert \geq
\vert F_j\cap F\vert $. Then $f(x)\not \in F$, otherwise
$R'_{\restriction F''_{j}}$, where $F''_j:=F'_j\setminus \{f(x)\}$, has type $\tau_i$ and $\vert F''_j\cap F\vert <\vert F_i\cap
F\vert$, contradicting the choice of $F_i$. Hence $f(x)\in
F'_j\setminus F$. Since $R'$ is $F$-monomorphic, the restriction
$R'_{\restriction F'_j \setminus \{f(x)\}}$ and $R'_{\restriction
F'_j\setminus \{x\}}$ are isomorphic. Since their types are respectively
$\tau_i$ and $\tau_j$, we have $i=j$.\hfill $\blacksquare$

\subsubsection{Directed Graphs}
Monomorphic decompositions can be easily described  in the case of directed graphs.

For that  we recall that a subset $A\subseteq V(G)$ of a a directed
graph is {\it autonomous} if for every $x,x'\in A, y\not \in A$, the
following two conditions holds:$$(x,y)\in E(G) \; \text{ if and only if}
\; (x',y)\in E(G)$$ $$(y,x)\in E(G) \; \text { if and only if } 
(y,x')\in E(G)$$.

The empty set, the one-element subsets of $V(G)$ and $V(G)$ itself are
autonomous. If there are no others, $G$ is {\it prime}.

We also recall that, if a directed graph $G$ is a lexicographical sum
$\sum_{i\in D}G_i$ of a family of directed graphs $G_i$, indexed by a
directed graph $D$, then provided that there are pairwise disjoint, the
$V(G_i)$ form a partition of $V(G)$ into automomous subsets. Conversely,
if the vertex set $V(G)$ of directed graph $G$ is partionned into
autonomous subsets, then $G$ is the lexicographical sum of the directed
graphs induced on the blocks of the partition.

\begin{theorem}\label{directedgraphs}
Let $G$ be a directed graph. Then $G$ has a finite monomorphic
decomposition if and only if $G$ is a finite lexicographical sum
$\sum_{i\in D}G_i$ of a family of directed graphs $G_i$, indexed by a
finite directed graph $D$, each $G_i$ being one of six kinds: an \emph{oriented acyclic
graph}, a \emph{clique}, an \emph{independent set}, a \emph{chain}, a \emph{reflexive clique}, or an
\emph{antichain}. Moreover, if $G$ decomposes into such a sum with all the
$V(G_i)$'s non-empty and if $D$ contains no $2$-element autonomous
subset, then both the monomorphic components of $G$ and the $V(G_i)$'s
wich are infinite coincide. 
\end{theorem}

\noindent{\it Proof.\enskip}
Let $W$ be a monomorphic component of $G$. If $W$ is infinite then
$G_{\restriction W}$ is of one of the six kinds mentionned in Theorem \ref{directedgraphs}.
Moreover, as it follows from implication $(iii)\Rightarrow (i)$ of
Theorem \ref{monomorphiccomp}, $W$ is autonomous. Hence, if $G$ has a finite
monomorphic decomposition, the partition of $V(G)$ into the infinite
blocks of the canonical decomposition of $G$ and of the singletons of
the remainder is a partition into autonomous sets, each of one of these
six kinds. If $G$ is a finite lexicographical sum $\sum_{i\in D}G_i$ of
a family of directed graphs $G_i$, each one  of these six kinds and if a monomorphic
component $W$ contains some $V(G_i)$ with $V(G_i)$ infinite then
$W=V(G_i)$. Otherwise, since $W$ is one of the six kinds above and $W$
is autonomous, the set $A:= \{j \in V(D): V_j\subseteq W\}$ is
autonomous and $D_{\restriction A}$ is of one of these six kind. Since
$D$ is finite, it contains a $2$-element autonomous subset. \hfill $\blacksquare$\vskip10pt

In the later case,  Lemma \ref {growthdegree} yields:

\begin{corollary}\label{alldegree}
If $G$ is a finite lexicographical sum $\sum_{i\in D}G_i$ of a family of
non-empty directed graphs $G_i$, indexed by a finite directed graph $D$,
each $G_i$ being either an oriented acyclic graph, a clique, an
independent set, a chain, a reflexive clique, or an antichain and if $D$
contains no $2$-element autonomous subset then $\varphi_{G}(n)\sim a
n^{k-1}$ where $k$ is the number of infinite $G_i's$.
\end{corollary}

A prime directed graph $D$ with $\vert V(D)\vert \geq 3$ cannot contain
a $2$-element autonomous subset. Since there are prime directed graphs
of arbitrarily large size, Corollary \ref {alldegree} yields examples of
directed graphs of polynomial growth of degree $k$ for every non
negative integer~$k$.

\subsubsection{Tournaments}\label{subsectiontournament} 
The monomorphic tournaments are the acyclic ones and the $3$-element
cycle. Hence from Theorem \ref{directedgraphs}, a tournament has a
finite monomorphic decomposition if and only if it is a lexicographical
sum of acyclic tournaments indexed by a finite tournament. We may
reformulate this in a simpler form.

An {\it acyclic component} of a tournament $T$ is a subset of $V(T)$
which is maximal w.r.t. inclusion among the acyclic autonomous subsets
of $V(T)$. Clearly, every acyclic autonomous subset is contained into an
acyclic component. As it is easy to see, the acyclic components of a
tournament form a partition of its vertex set. It follows that a
tournament is a lexicographical sum of acyclic tournaments indexed by a
finite tournament if and only if it has only finitely many acyclic
components.

With Y.Boudabbous \cite{boudabbous}, we identified twelve infinite tournaments and proved that an infinite tournament $T$ is a finite lexicographical sum of acyclic tournaments if and only if  $T$ embeds none of these tournaments. The growth of the profile of each of these tournaments being exponential, we deduced from Theorem \ref{theorem.quasipolynomial}  and Lemma \ref{growthdegree} the following dichotomy result.   

\begin{theorem} \label{thetheorem} \cite{boudabbous} The growth of the profile of a tournament $T$ is either  polynomial, in which  case $T$ is a lexicographic sum of acyclic tournaments indexed by a finite tournament,  or it is  at least exponential.
\end{theorem}




There are prime tournaments of arbitrarily large finite size, hence,
according to Corollary \ref{alldegree} there are tournaments of
arbitrarily large polynomial growth. We give below some examples of
small growth.

\begin{examples}\label{exampletournament}
\begin{enumerate}
\item If $T$ is an acyclic tournament, then $\varphi_T(n)=1$ for every integer $n$, $n\leq \vert T\vert$. Conversely, if $\vert T\vert \not =3$ and $\varphi_T(3)=1$ then $T$ is acyclic.

Let $\omega$ be  the tournament made of the integers, with the strict ordering, that is  $\omega:=(\N, \{(p,q)\in \N^{2} : p < q\})$, and let $\omega^*$ be its dual. Note that according to the theorem of Ramsey, every  infinite tournament contains a subtournament which is isomorphic to $\omega$ or to $\omega^*$.

Let  $C_3:= (\{0,1,2\}, \{(0,1),(1,2), (2,0)\})$ be the $3$-element cycle. For $i:=1,2, 3$, let
$T_i$ be the tournament obtained by replacing $i$ vertices of $C_3$  by $\omega$.

\item   $\varphi_{T_1}(0)=\varphi_ {T_1}(1)=\varphi_{T_1}(2)=1$, $\varphi_{T_1} (n)= 2$ for all $n\geq 3$.

\item   $\varphi_{T_2}(0)=\varphi_ {T_2}(1)=\varphi_{T_2}(2)=1$, $\varphi_{T_2} (3)=2$,  $\varphi_{T_2} (n)= n-2$ for all $n\geq 4$ and $H_{\varphi _{T_{2}}}= \frac{1-x+x^3-x^4+x^5}{(1-x)^2}$.

\item  $T_{3}=\omega.C_3$ and  $H_{\varphi_{T_{3}}}=\frac{1-x^2+x^5+x^6}{(1-x)(1-x^2)(1-x^3)}=1+x+x^2+2x^3+ 2x^4+3x^5+5x^6+6x^7+8x^8+11x^9+13x^{10}+16x^{11}+ O(x^{12})$.

\end{enumerate}

\end{examples}

The profile of tournaments goes largely beyond polynomials. For an example of exponential profile, let
$C_3. \omega$ be the lexicographical sum of infinitely many copies of
$C_3$. Then the generating serie of the profile is
$H_{\varphi_{C_{3.\omega}}}=\frac{1}{1-x-x^3}$ hence the profile is
exponential (see sequence A000930 in Encyclopedia of integer
sequences (URL: http://www.research.att.com).).


\section{The Age Algebra of Cameron}
P.~J.~Cameron~\cite{cameronbook,
cameron.1997} associates to $\age(R)$,  the age of  a relational structure $R$,  its \emph{age
  algebra}, a graded commutative algebra $\K.\age(R)$ over a field
$\K$ of characteristic zero.  He shows that if  $\varphi_R$ takes only finite values, then the dimension of $\K.\age(R)_n$, the
homogeneous component of degree $n$ of $\K.\age(R)$,  is
$\profile_R(n)$.  Hence, in this case,  the generating series of the profile  is simply the Hilbert series of $\K.\age(R)$.

P.J Cameron mentions several interesting examples of algebras which turn
to be age algebras. The most basic one is the {\it shuffle algebra} on
the set $A^*$ of words on a finite alphabet $A$ \cite{lothaire}. Indeed,
as mentionned at the end of Subsection 2.3, $A^*$ is the age of the relational
structure $(\Q , (U_a)_{a\in A})$ where the $U_a$'s are unary relations
forming a coloring of $\Q$ into distinct colors, in such a way that
between two distinct rational numbers, all colors appear. And the
shuffle algebra is isomorphic to the age algebra of $(\Q , (U_a)_{a\in
A})$.

There are several reasons to associate a graded algebra to an age.
We will examine some in the next subsections. For the ease of our
discussion, we recall the presentation of the age algebra via the set
algebra (see \cite{cameron3}).

\subsection{The Set Algebra} 

Let $E$ be a set and let $[E]^{<\omega}$ be the set of finite subsets of $E$
(including the empty set $\emptyset$). Let $\K$ be a field and
$\K^{[E]^{<\omega}}$ be the set of maps $f:[E]^{<\omega}\rightarrow \K$. Endowed with the usual
addition and scalar multiplication of maps, this set is a vector space over $\K$. Let
$f,g\in
\K^{[E]^{<\omega}}$ and $Q\in [E]^{<\omega}$. Set
\begin{equation}
f g(Q) = \sum_{P \in [Q]^{<\omega}} f(P)g(Q\setminus P)
\end{equation}.
With this operation added, the above vector space becomes a  $\K$-algebra. This algebra is commutative and it has a
unit, denoted by $1$. This is the map taking the value $1$ on the empty set and the value $0$
everywhere else. The \emph{set algebra} is the subalgebra
 made of the
maps such that $f(P)=0$ for every $P\in [E]^{<\omega}$ with $\vert
P\vert $ large enough. This algebra is graded, the homogeneous component
of degree $n$ being made of maps which take the value $0$ on every
subset of size different from $n$.

 Let $\equiv$ be an equivalence relation on $[E]^{<\omega}$. A map $f:[E]^{<\omega}\rightarrow \K$ is
{\it $\equiv$-invariant},  or  briefly, {\it invariant},  if $f$ is constant on each equivalence class.
Invariant maps form a subspace of the vector space $\K^{[E]^{<\omega}}$.
We give a condition which insure that they form a subalgebra too.

An equivalence relation on $[E]^{<\omega}$ is {\it hereditary} if every pair $D,D'$ of equivalent
elements  satisfies the following conditions:
\begin{enumerate}\item $\vert D\vert= \vert D'\vert=d$ for some $d$.
\item $\vert \{X\subseteq D:X\equiv B\}\vert =
\vert \{X\subseteq D':
X\equiv B\}\vert \; \mbox{ for every}\;  B\subseteq E$
\end{enumerate}

Hereditary equivalence where introduced in \cite {pouzetrosenberg}.

With N.Thi\'ery, we proved:
 \begin{proposition}
Let $\equiv$ be an hereditary equivalence relation on $[E]^{<\omega}$.  Then the
product of two invariant maps is invariant. Thus  the set of invariant maps
$f:[E]^{<\omega}\rightarrow \K$ form a subalgebra of $\K^{[E]^{<\omega}}$.
\end{proposition}

Let $R$ be a relational structure with domain $E$.  Set $F\equiv_R F'$
for $F,F'\in [E]^{<\omega}$ if the restrictions $R\restriction_{F}$
and $R\restriction_{F'}$ are isomorphic. The resulting equivalence on
$[E]^{<\omega}$ is hereditary. Let $\K.\age(R)$ be the intersection of the subalgebra of
$\K^{[E]^{<\omega}}$ made of invariant maps with the set algebra. This is a graded  algebra, the \emph{age
  algebra} of Cameron.

  The name comes from the fact that this algebra depends only upon the age of $R$. Indeed,
if $R'$ is a relational structure with domain $E'$ such that $\age(R')=
\age(R)$ then $\K.\age(R')$ identifies to $\K.\age(R)$. To see that,
associates first to every indicator function $\chi_{S}$ of an
equivalence class $S$ for $\equiv_R$ the indicator function $\chi_{S'}$ of
the equivalence class $S'$ for
$\equiv_{R'}$ such that $R_{\restriction P}$ is isomorphic to
$R'_{\restriction P'}$ for some $P\in S$ and some $P'\in S'$. Next
associate to every linear combination (finite or not) of indicator
function the linear combination, with the same coefficients, of their
images.

If $\varphi_R$ takes only integer values, $\K.\age(R)$ identifies with
the set of (finite) linear combinations of members of $\age(R)$. This
explain the fact that, in this case, $\varphi_R(n)$ is the dimension of
the homogeneous component of degree $n$ of $\K.\age(R)$.
In a special case, we have

\begin{theorem}\label{polyalgebra2}\cite {pouzetthiery}
If $R$ has a monomorphic decomposition into finitely many blocks $E_1,
\dots, E_k$, all infinite, then the age algebra $\K.\age(R)$ is a
polynomial algebra, isomorphic to a subalgebra $\K[x_1, \dots, x_k]^ R$
of $\K[x_1, \dots, x_k]$, the algebra of polynomials in the
indeterminates $x_1, \dots, x_k$.
\end{theorem}

\subsection{Behavior of the Profile}\label{behavior}

In the frame of its age algebra,  Cameron gave the following proof of the fact  that the profile does not decrease.

Let $R$ be a relational structure on a set $E$,   let $e:=\sum_{P\in [E]^{1}}P$ (that we could view as the sum of isomorphic types of the one-element restrictions of $R$) and $U$ be the subalgebra generated by $e$.  Members of $U$ are of the form $\lambda_me^m+\cdots+\lambda_1e+\lambda_01$ where $1$ is the isomorphic type of the empty relational structure and $\lambda_m,\dots, \lambda_0$ are in $\K$. Hence $U$ is graded, with $U_n$, the homogeneous component of degree $n$, equals to $\K.e^n$.

 \begin{theorem}\label{cameron}
If $R$ is infinite then for every $u\in \K.\age(R)$, $eu=0$  if and only if $u=0$
 \end{theorem}

This innocent looking result implies that $\varphi_R$ is non decreasing.
Indeed, the image of a basis of $\K.\age(R)_n$ by multiplication by
$e^m$ is an independent subset of $\K.\age(R)_{n+m}$.

\subsection{\bf Finite generation}

If  a graded
 algebra $A$ is finitely generated, then,  since $A$ is a quotient of the polynomial ring $\K[x_{1},\dots, x_{k}]$,  its Hilbert function is bounded above by a polynomial. In fact, as it is well known, its Hilbert  series is a fraction of  form $\frac{P(x)}
{(1-x)^d}$, thus of the form given in (\ref{hilbertpoly}) of subsection \ref{subsectionpoly}. Moreover, one can choose a numerator  with  non-negative coefficients whenever the algebra is Cohen-Macaulay.  Due to Problem \ref{growthseries}, one could be tempted to conjecture that  these  sufficient   conditions are necessary in the case of age agebras.
Indeed, from Theorem \ref{cameron} one deduces easily:
 \begin{theorem}\label{finitemodule} 
The profile of $R$ is bounded if and only if $\K.\age(R)$ is finitely
generated as a module over $U$. In particular, if one of these
equivalent conditions holds, then $\K.\age (R)$ is finitely generated
 \end{theorem}

 But this case is exceptional. Indeed, on one hand,
as we have mentionned in \ref{subsectiontournament}, there are
tournaments whose profile has arbitrarily large polynomial growth rate.
On an other hand, with N.Thiery we proved:
 \begin{theorem} 
The age algebra of a tournament  is finitely generated if and only if the profile  is bounded.
 \end{theorem}

The argument is simple and illustrates the  above notions.
Let $T$ be  a tournament. Suppose that  $\varphi_T$ is  bounded, then by Theorem \ref{finitemodule},  $\K. \age(T)$ is finitely generated.  Conversely, suppose that $\K. \age(T)$ is finitely generated. Then,  as mentionned in the introduction of this subsection, $\varphi_T$ is bounded above by a polynomial.  Apply Theorem \ref{thetheorem}. Necessarily,  $T$ is a  lexicographical sum of acyclic tournaments indexed by a finite tournament, a fact that we may write  $T= \sum_{i\in D}A_i$, where each $A_i$ is an acyclic tournament and $D$ is a finite tournament. If $\varphi_T$ is  not bounded by a constant then, as one can easily see, $D$ contains   a $3$-element subset $A:= \{i,j,k\}$ such $D_{\restriction A} $ is a cycle and  $A_j$ and $A_k$ are infinite. Pick an element $a_i\in A_i$, set $E':= \{a_i\}\cup A_j\cup A_k$  and set $T':= T_{\restriction E'}$. Since $T'$ is an induced substructure of $T$,  the algebra $\K. \age(T')$ is a quotient of $\K. \age (T)$. Thus, with our hyp
 othesis that $\K. \age (T)$ is finitely generated, $\K. \age (T')$ must be finitely generated too.
Let us see that this is not the case. Indeed, notice that  if $T'$ is an arbitrary infinite tournament then  in the age algebra $\K.\age (T')$ we have $e^n=n!(a_n+b_n)$
where  $a_n$ is the isomorphic type of an acyclic tournament on $n$ vertices and $b_n$ is the sum of isomorphic types of induced subtournaments of $T'$ on $n$ vertices which contain a triangle.
It follows then that every element $x\in \K.\age(T')$ is a sum $x:=a(x)+u(x)$ where $a(x)$ is a linear combination of types, each containing a cycle, and $u(x)\in U$. Thus, if  $g_1, \dots, g_k, 1$ generate $\K.\age(T')$ then $a(g_1), \dots, a(g_k), e$ generate it also. We apply this fact to our  tournament $T'$.  Since $T'$ does not contain two disjoint cycles, we have  $a(g_i)a(g_j) =0$ for every $1\leq i,j\leq k$. Thus, $a(g_1), \dots, a(g_k),1$ generate $\K. \age(T)$ as a module over $U$. According to Theorem \ref{finitemodule},  $\varphi_{T'}(n)$ is bounded.  But, $\age(T')=\age (T_2)$ where $T_2$ is the  example given in Item 3 of Examples \ref{exampletournament}. Hence according to the computation given there, $\varphi_{T'}(n)=n-2$ for $n\geq 4$. This gives a contradiction. Thus  $\varphi_T$ must be bounded by a constant, as claimed.

\subsection{The Behavior of the Hilbert Function; a Conjecture of Cameron}

 Cameron \cite{cameron4} made an important observation about the behavior of the Hilbert fonction.

\begin{theorem}\label{hilbertnonzero} 
Let  $A$ be a graded algebra over  an algebraically closed field of
characteristic zero.  If
$A$ is an integral domain  the values of the Hilbert function $h_A$ satisfy the inequality
\begin{equation}\label{eqhilbert}
h_A(n)+h_A(m) -1\leq h_A(n+m)
\end{equation}
for all non-negative integers $n$ and $m$. 
\end{theorem}

In 1981, he made the
conjecture that if $R$ codes  a permutation groups with no finite
orbits then the age algebra over $C$ is an integral domain \cite{cameron81}
(see also \cite {cameron3} p.~86). I solved it positively in 2002 in a
slightly more general setting:

\begin{theorem}\label{agealgebranozero} 
Let $R$ be a relational structure with possibly infinitely many non
isomorphic types of $n$-element substructures. If the kernel of $R$ is
empty, then $\K.\age(R)$ is an integral domain.\end{theorem}

Since the kernel of a relational structure $R$ encoding a permutation
group $G$ is the union of its finite orbits, if $G$ has no finite orbit,
the kernel of $R$ is empty. Thus from this result, $\K.\age(R)$ is an integral domain, as conjectured by Cameron.

At the core of the solution is this lemma:
\begin{lemma}\label{main}
Let $m,n$ be two non negative integers. There is an integer $t$ such that for every set $E$, every field $\K$ with characteristic zero, every pair of maps 
$f:[E]^m\rightarrow
\K$,
$g:[E]^n\rightarrow
\K$, if $f g(Q) := \sum_{P \in [Q]^{m}} f(P)g(Q\setminus P)=0$ for every $Q\in [E]^{m+n}$
 but $f$ and $g$ are not identically zero, then $f$ and $g$ are zero on
$[E\setminus S]^m$ and
$[E\setminus S]^n$ where
 $S$ is a finite subset of $E$  with size at most $t$ \end{lemma}

If the age is inexhaustible, then in order to
prove that there is no zero divisor, the only part of the lemma we need to apply is the assertion
that $S$ is finite.

The fact that
$S$ can be bounded independently of $f$ and
$g$, and the value of the least upper bound  $\tau(n,m)$,  seem  to be of independent
interest.  The only exact value we know is
$\tau(1,n)=2n$, a fact which amounts to a weighted version of Theorem \ref{kantor}.
Our existence proof of $\tau(m,n)$ yields  astronomical upper bounds. For example,  it gives $\tau (2,2)\leq 4R_k^2(4)+1$, where $k:=5^{56}$ and $R_k^2(4)$ is the Ramsey  number equals to the least integer $p$ such that for every colouring  of the  pairs of $\{1, \dots, p\}$ into $k$ colors there are four integers whose all pairs have the same colour.  The only lower bound we have is  $\tau(2,2)\geq 7$ and more generally $\tau(m,n)\geq (m+1)(n+1)-2$.
 We cannot preclude a
extremely simple upper bound for $\tau(m,n)$, eg quadratic in
$n+m$.

For example, our 1971 proof of Theorem \ref{increaseinfinite} consisted to show  that
$\varphi_{R}(n)\leq
\varphi_{R}(n+1)$ provided that $E$ is large enough, the size of $E$ being bounded by some Ramsey number, whereas,
according to Theorem \ref{kantor},
$2n+1$ suffices \cite{pouzet76}.

Our proof of Lemma \ref{main} is built on the integrity of the age
algebra of infinite multichainable relational structures and Ramsey's
theorem. Here is a sketch (details will appear in \cite{pouzet02}).
 We reduce the first part to the integrity of a shuffle algebra. For
this, let $R$ be a multichainable relational structure and let $V$ and
$L$, with $L$ infinite, as in \ref{subsectionmultichain}. Consider the vector
space over $\K$ spanned by the words whose letters are non-empty
subsets of a finite set $V$, the shuffle $u \shuffle v$ of two words $u$ and $v$ being the sum of
all words $w$ which are the disjoint union of one occurrence of $u$ and one occurrence of $v$ (eg
if
$V:=\{a,b\}$ then $\{a\}\shuffle
\{b\}=\{a\}\{b\}+\{b\}\{a\} +\{a,b\}$).  Using  a lexicographical order, one  can show that the resulting algebra is an
integral domain \cite{radford}. This algebra is the age algebra of a
relational structure $M$ of a special form: the product of a relational
structure on $V$ by an infinite chain $L$ (see subsection \ref
{subsectionproduct}). Since the isomorphic classes of $R$ contain those
of $M$, the age algebra of
$R$ is an integral domain.

Next, we consider a
pair $f$, $g$ of maps as in Lemma \ref{main}, we suppose $m\leq n$ and that  the subsets $S$ for which the conclusion of the lemma holds are "large" (that is $\vert S\vert \geq \mu (m,n)$  where $\mu(m, n)$ is a number depending on $\tau(m,n-1)$, $\tau(m-1, n)$, and  $r(n+m):=R^{n}_{k}(m+n)$, where $k:=5^s$ and $s:= {(m+n)n \choose
m}+{(m+n)n\choose n}$. We observe then that there are at least
$r(m+n)$ pairwise disjoint $m+n$-element subsets $V_{0},\dots,V_{i},\dots, V_{r(m+n)-1}$ of $E$, each  satisfying
$f(P_{i})g(V_{i}\setminus P_{i})\not =0$ for some $P_{i}\in [V_{i}]^m$. Supposing that the union $E'$ of these subsets is
the product of an $m+n$ element set $V$ by a $r(m+n)$-element chain $L'$, we encode the action of $f$ and $g$ on $E'$ by a
relational structure $R'$ made of a $m$-uniform hypergraph and a $n$-uniform hypergraph whose hyperedges are coloured in at
most four
colors. Dividing the
$m+n$-element subsets of $L'$ into equivalence classes, Ramsey's theorem provides a $m+n$-element subchain
$L''$ such that the isomorphic classes of $R'_{\vert V\times L''}$ contains those of the
product of $V$ by $L''$. From the fact that the shuffle algebra built on $V$ is an integral domain we deduce that either the
restriction of $f$ to $[V\times L'']^m$ or the
restriction of $g$ to $[V\times L'']^n$  vanishes. But this is not the case. Hence
$S$ cannot be large.

\subsection{Initial  segments of an age and ideals of a ring}

As stated in Theorem \ref{profilepouzet2}, if the profile of a
relational structure $R$, with bounded signature or finite kernel, has a
polynomial growth, then it is almost multichainable. As we will see in
Theorem \ref{productwqo}, if a relational structure $R$ is almost
multichainable, its age $\age(R)$, ordered by embeddability, is {\it
well-quasi-ordered}, that is every final segment of $\age (R)$ is
finitely generated, which amounts to the fact that the collection
$F(\age(R))$ of final segments of $\age(R)$ is noetherian, w.r.t. the
inclusion order (see subsection \ref{subsectionwqo}). Final segments
play for posets the same role than ideals for rings. Noticing that an
age algebra is finitely generated if and only if it is noetherian, we
are lead to have a closer look at the relationship between the basic
objects of the theory of relations and of ring theory, particularly ages
and ideals.

We mention the following result which will be incorporated into a joint paper with  N.Thi\'ery.
\begin{proposition} 
Let $\mathcal A$ be the age of a relational structure $R$ such that the
profile of $R$ takes only finite values and $\K.\age $ be its age
algebra. If $\mathcal A'$ is an initial segment of $\mathcal A$ then:
 \begin{enumerate} 
 \item[{(i)}] The vector subspace $J:= \K.(\mathcal A\setminus \mathcal A')$ spanned by $\mathcal A\setminus \mathcal A'$ is an ideal of $\K.\age $. Moreover, the quotient of $\K. \age$ by $J$ is  a ring isomorphic to  the ring $\K.\mathcal A'$.
 \item[(ii)] If this ideal is irreducible then $\mathcal A'$ is a subage of $\mathcal A$.
  \item[(iii)] This is a prime ideal if and only if $\mathcal A'$ is an inexhaustible age.
 \end{enumerate}
\end{proposition}
The proof of Item (i) and Item (ii) are immediate. The proof of Item
(iii) is essentially based on Theorem \ref{agealgebranozero}.

According to Item (i), $F(\mathcal A)$ embeds into the collection of ideals of $\K.\mathcal A)$. Consequently:
\begin{corollary} 
If an age algebra is finitely generated then the age is well-quasi-ordered by embeddability.
\end {corollary}

\begin{problem}
How the finite generation of an age algebra translates in terms of embeddability between members of the age?
\end{problem}

\section{Tools for Classifying Profiles}

Beside graded algebras, Hilbert series and rudimentary algebraic
geometry, the tools we use to classify profiles can be divided into the
following categories:
\begin{enumerate}
\item[--] Orders, well-founded orders, well-quasi-orders;
\item[--] Ramsey Theorem;
\item[--]  Compactness theorem;
\item[--] Combinatorial properties of the kernel.
\end{enumerate}
In this section, we show the role of these tools in the proof of Theorem
\ref{profilpouzet1} and Theorem \ref {profilepouzet2}. For reader
convenience, we record in the first subsection below some notions and
notations we are using in this paper.

\subsection{Basic Notions, Relational Structures,  Embeddability  and Ages}

We use standard set-theoretical notations. If $E$ is a set, $\vert
E\vert$ denotes its cardinality. If $n$ is an integer, $[E]^n$ denotes
the set of $n$-element subsets of $E$; whereas $E^n$ denotes the set of
$n$-tuples of elements of $E$. An {\it $n$-ary relation on $E$} is any
subset $\rho$ of $E^n$. As said in subsection \ref{defexamples}, a pair
$R:= (E,
(\rho_i)_{i \in I})$ made of a set $E$ and of a family of $m_i$-ary
relations $\rho_i$ on $E$ is a {\it relational structure of signature} $\mu:=(m_i)_{i\in I}$. If $I'$ is a subset of $I$, the relational structure $R^{I'}:= (E,
(\rho_i)_{i \in I'})$ is a {\it reduct} of $R$.

Let $R:= (E, (\rho_i)_{i\in I})$ and $R':= (E', (\rho'_i)_{i\in I})$ be two
relational structures of the same signature $\mu:=(m_i)_{i\in I}$.

 A map $f: E\rightarrow E'$ is \emph{an isomorphism  from} $R$ {onto} $R'$ if
 \begin{enumerate}
 \item $f$ is bijective,
 \item $(x_1, \dots, x_{m_i})\in \rho_{i}$ if and only if $(f(x_1), \dots, f(x_{m_i}))\in \rho'_{i}$ for every $(x_1, \dots, x_{m_i})\in E^{m_i}$, $i\in I$.
 \end{enumerate}
A map $f$ from a subset $A$ of $E$ onto a subset $A'$ of $E'$ is a
\emph{local isomorphism} or a \emph{local embedding} of $R$ into $R'$ if
$f$ is an isomorphism from $R_{\restriction A}$ onto $R'_{\restriction
A'}$. If $A= E$, this is \emph{an embedding from $R$ into $R'$}.

 We say that  $R$ is \emph{isomorphic}, resp. \textit{embeddable}\index{embeddable},  to $R'$, resp. into $R'$, and we write $R\cong R'$, resp.
$R\leq R'$, if there is an isomorphism, resp. an embedding,  from $R$ onto $R'$, resp. into $R'$. Clearly, $R$ is embeddable into $R'$ if and only if $R$ is isomorphic to some restriction of $R'$.
To each relational structure $R$ one may associate an \emph{isomorphic
type $\tau(R)$}, in such a way that $ R\cong R'$ if and only if
$\tau(R)=\tau(R')$. The collection of isomorphic types of finite
relational structures of signature $\mu$ is a set, that we denote
$\Omega_{\mu}$.
The embeddability relation is a  quasi-order on the
class of relational structures of signature $\mu$. It induces a
quasi-order on the class of their isomorphic types. On the
set $\Omega_{\mu}$ this is an ordering.

The \textit{age} of a relational structure  $R$
is the
set ${\mathcal A}(R)$
of isomorphic types of the  restrictions of $R$ to the finite subsets of its domain.

Ages  can be (almost) characterized in terms of ordered sets (posets). Let us recall that if $P$ is a poset,
an {\it initial segment} of $P$ is a subset $I$ such that $x\in P$, $y\in I$ and $x\leq y$ imply $x\in
I$. An {\it ideal} is a non-empty initial segment which is {\it up-directed}, that is $x, y \in
I$ imply $x, y\leq z$ for some $z\in I$.

The following characterization is essentially  due to R.
Fra\"{\i}ss\'e 1954 (see \cite{fraissetr}).

\begin{lemma}\label{121} $a)$ Let $R$ be a relational structure of
  arity $\mu$ and size $\kappa$, then ${\mathcal A}(R)$, the age of $R$, is an
ideal of
  $\Omega_\mu$, whose size is at most $\kappa$ if $\kappa$ is infinite
  and $2^{\kappa}$ otherwise;\\
$b)$ Conversely, let ${\mathcal C}$ be an ideal of $\Omega_\mu$; if
${\mathcal C}$ is countable, then ${\mathcal C}$ is the of age a
countable $\mu$-ary relational
structure.\\
 In particular:\\
$c)$ If $\mu$ is finite, then a subset ${\mathcal C}$ of
$\Omega_{\mu}$  is the age of a  $\mu$-ary relational structure if
and only if ${\mathcal C}$ is an ideal of $\Omega_\mu$.\end{lemma}

The equivalence stated in c) does not hold without some condition on
$\mu$. An example was obtained with
C.~Delhomm\'e and mentionned in  \cite{pouzetsobranisa}, see \cite{DPSS} for other examples and further developments.

Still, $\Omega_{\mu}$ is an age -and the largest one. As a poset, it decomposes into
{\it levels}, the levels of a poset $P$ being defined inductively by the formula
$P_{n}:= Min (P\setminus \cup \{ P_{n'}: n'<n \})$. The $n$-th level is the set of isomorphic
types  of relational structures on $n$ elements. The same holds for the  age of a
relational structure
$R$, hence the profile $\varphi_R (n)$ is simply the function which assign to each $n$ the
number of elements of the $n$-th level of ${\mathcal A}(R)$  (in the general theory of
posets, the numbers of elements of the  levels of a poset $P$ are the {\it Whitney numbers}- of
the second kind- of $P$). In the study of their  growth, we have only to consider profiles taking
finite values. Hence, according to the above characterization,  our study is about ideals of
$\Omega_{\mu}$, each one  with finite levels.  A  characterization  in order-theoretical terms,
of  these  ideals eludes us,   the study of the profile
is just a bare approach.

\subsection{The Height Function on a Well-Founded Poset,  
the Height of\newline an age and its Relation with the Profile}

Let $P$ be a poset;  $P$ is \emph{well-founded} if every non-empty subset $A$ of $P$ contains some minimal element $a$ (that is there is no $b\in A$ with $b<a$).  The \emph{height function} on a poset $P$ associates to an element  $x\in P$ an ordinal number  $h(x,P)$ is such a way that:
$$h(x,P):= Sup \{h(y, P)+1: y<x\}$$

The ordinal $h(x, P)$ is   the \emph{ height of $x$ in $P$}.

With this definition, $h(x, P)=0$ if and only if $x$ is minimal in $P$;
also $h(x, P)$ is defined if and only if $\downarrow x: =\{y\in P: y\leq
x\}$ is well-founded.

Let
$\mathcal {A}$ be an ideal of $\Omega_\mu$ and  let  $\mathcal {J}(\mathcal {A})$ be the set of ideals included into
$\mathcal {A}$. The  {\it height} of $\mathcal {A}$,
denoted by
$H({\mathcal A})$,  is the {\it height} of $\mathcal {A}$ in $\mathcal {J}(\mathcal {A})$.
With this definition, $H(\mathcal A) = n$ with $
n < \omega$, if and only if  $\mathcal A=\mathcal A (R)$ for some  $R$ on $n$ elements.
Furthermore, $H(A) =
\omega$ if and only if $\mathcal A$ is infinite and every ideal properly included
into
${\mathcal A}$ is finite.

The characterization of  ideals of height $\omega$ is given by the following result:

\begin{theorem}[\cite{pouzetrpe}]
Let $\mathcal A$ be an ideal  of $\Omega_{\mu}$.  Then the following properties are equivalent:
\begin{enumerate}
\item[{(i)}] $H(\mathcal A)= \omega$.
\item[(ii)] $\mathcal A$ is the age of an infinite monomorphic relational structure $R$ (that is $\varphi_R(n)=1$ for all $n$).
\item[(iii)] $\mathcal A$ is the age of an infinite chainable relational structure.
\end{enumerate}
\end{theorem}
The essential part is the implication $(i)\Rightarrow (iii)$. For that,  observe that $\mathcal A$ is countable.  Hence, from Lemma \ref{121} b), $\mathcal A$ is the age of some countable relational structure $R$. From Theorem \ref{chainablerestriction} and Compactness theorem of first order logic, $R$ is chainable.

This generalizes:
\begin {theorem}
[\cite{pouzetrpe}]
Let $\mathcal A$ be an ideal  of $\Omega_{\mu}$.  Then the following properties are equivalent:
\begin{enumerate}
\item[{(i)}]  $H(\mathcal A ) = \omega + p$,  for some $p<\omega$.
\item[(ii)]  $\mathcal A$ is the age of an infinite  relational structure with a bounded profile.
\item[(iii)] $\mathcal A$ is the age of an infinite  almost chainable relational structure.
\end{enumerate}
\end{theorem}

 This result opens the road for proving Theorem \ref{profilpouzet1}.
 Indeed, as it will appear at  the end of this section:
\begin{theorem}\label{heightprofile}
If $R$ has a finite kernel then the growth of
$\varphi_{R}$ is polynomial of degree $k$ if and only if  $\mathcal A(R)$ has an height and
$H(\mathcal A(R)) = \omega(k+1) + p$ for some integer $p$.
\end{theorem}
At this point, we are still far away of a proof. More tools are needed.

\subsection{Well-Quasi-Ordering, Higman Theorems,\newline
 Very-Well-Quasi-Ordered Classes}\label{subsectionwqo}

A poset $P$ is {\it well-quasi-ordered}, in brief w.q.o., if every
non-empty subset contains finitely many minimal elements(this number
being non-zero). A final segment $F$ of a poset  $P$ is {\it finitely generated}  if for some finite subset $K$ of $P$, $F$ equals the set $\uparrow K:=\{y\in P: x\leq y\;  \text{for  some} \; x\in K\}$.

A basic result due to Higman \cite{higman52} is the following.

\begin{theorem}\label {higman1} For a poset  $P$, the following properties are equivalent:
\begin{enumerate}
 \item $P$ is  w.q.o.
 \item  $P$ contains  no infinite strictly descending chain and no infinite antichain.
\item  Every final segment of $P$ is finitely generated.
\item The set $I(P)$ ordered by inclusion and made  of initial segments of $P$ is well-founded.
\end{enumerate}
\end{theorem}

\begin{corollary}
If an ideal $\mathcal A$ of $\Omega_{\mu}$ is w.q.o. then $\mathcal J(\mathcal A)$ is well-founded, hence $H(\mathcal A)$ is defined.
\end{corollary}
Non-w.q.o posets for which the collection of ideals is well-founded  abund. Not a single example of non-w.q.o age whose the collection of subages is well-founded as been discovered yet.

With Zorn lemma, one may observe that  if a poset $P$ is not w.q.o.  it contains some  final segment $F$ which is not finitely generated and maximal w.r.t. inclusion. If follows that its complement  $P\setminus F$ is w.q.o.

This fact applies nicely to $\Omega_{\mu}$. For that,  let us mention that a \emph{bound} of  a relational structure $R$ is any  relational structure $S$ on a finite set $F$ such that $S$ does not embed into $R$ but $S_{-x}:= S_{\restriction F\setminus \{x\}}$ embeds into $R$ for every $x\in F$. A {\it bound} of an initial segment $\mathcal C$ of $\Omega_{\mu}$ is any minimal element of  $\Omega_{\mu}\setminus \mathcal C$. Clearly, the bounds of $\mathcal A(R)$ are the isomorphic types of the bounds of $R$.
\begin{lemma}\label {nash-williams} 
If an initial segment $\mathcal C$ of $\Omega_{\mu}$ is level finite and is not w.q.o., then it contains an ideal which has infinitely many bounds in $\mathcal C$.
\end{lemma}

Let us recall that a {\it word} is a finite sequence of {\it letters}
and that a word $u$ is a {\it subword } of a word $w$ if $u$ can be
obtained from $w$ by erasing some letters of $w$. We have the following
theorem of Higman \cite{higman52}.

\begin{theorem}\label{higman2}
The set $A^{\ast}$ made of words on a finite alphabet $A$ is w.q.o. with
respect to the subword ordering. 
\end{theorem}

A class ${\mathcal C}$  of relational structures is  {\it very-well-quasi-ordered}, in short v.w.q.o., if for every
integer
$m$ the class ${\mathcal C}_{m}$ made of  $R\in {\mathcal C}$ added of  $m$
unary relations is w.q.o for the embeddability relation
We will need  the following result.
\begin{theorem}\label {pouzet 73}
Let  $\mathcal C$ be a class  of finite relational structures,
then:
\begin{enumerate}
\item  $\mathcal C$ is  v.w.q.o iff $\downarrow \mathcal C$, its downward closure, is
v.w.q.o.
\item If $\downarrow C$ is v.w.q.o. then all the ages it contains are almost inexhaustible.
\item If $\downarrow
\mathcal C$ is v.w.q.o.  and all of its members  have  the same finite
signature
$\mu$  then it has only finitely many bounds \cite{pouzet 73}.
\end{enumerate}
\end{theorem}

An other important result on w.q.o.  is  this (see  \cite{wolk67} for the first part, \cite {dejongh-parikh} for the second).

\begin{theorem} 
If a poset $P$ is w.q.o. then  all the linear extensions of $P$
are well-ordered and there is one having the largest possible order type.
\end{theorem}

 This largest order type, denoted $o(P)$,  is  the
{\it ordinal length} of $P$.

\begin{lemma}[\cite{dejongh-parikh}]\label {lengthalphabet}  
If $A$ is an alphabet made of  $k$ letters then $o(A^*)= \omega^{\omega^{k-1}}$.
\end{lemma}

The computation for ages (see \cite{pouzetsobraniia},  \cite{sobranietat}) yields:

\begin{theorem}
\label{longueur}If ${\mathcal J}({\mathcal A}(R))$ is w.q.o then $o({\mathcal A}(R))=\omega^\alpha\cdot q$
where $\alpha$  is such that   $\omega
\cdot\alpha\leq H({\mathcal A}(R))<\omega\cdot (\alpha +1)$ and  $q$ is the number
of ages included into ${\mathcal A}(R)$ whose height is between  $\omega
\cdot\alpha$ and $\omega
\cdot (\alpha +1)$.
\end{theorem}

\subsection{Kernel, Almost Inexhaustibility, Height and Profile}
Most of the properties of the kernel (defined in subsection \ref{jump})are based on the following simple lemma (see \cite{pouzetrm} for finite signature and 3 of Lemma 2.12 of \cite{pouzetsobranisa} for the general case).
\begin{lemma} \label{minusab}
$\mathcal {A}(R_{\restriction E\setminus \{a\}})=\age (R)\Longrightarrow \mathcal {A}(R_{\restriction E\setminus \{a, b\}})=\mathcal {A}(R_{\restriction E\setminus \{b\}})$ for all $a,b\in E$
\end{lemma}

From this lemma, one immediately gets the following:

\begin{corollary}
A relational structure $R$ has an empty kernel if and only if its age $\age(R)$ has the disjoint embedding property.
\end{corollary}

Instead of saying that an age has the   disjoint embedding property, we say that it is {\it inexhaustible}.    We say that  a relational structure $R$ is  {\it age-inexhaustible} if $K(R)$ is empty and
{\it almost age-inexhaustible } if  $K(R)$ is finite.
\begin{lemma}\label{limit}
Let $\mathcal A$ be an infinite inexhaustible age. Then:
\begin{enumerate}
\item  For every age $\mathcal A'$ included into $\mathcal A$ there is an infinite strictly  increasing sequence of  ages included into $\mathcal A$ such that $$\mathcal A'=\mathcal A_0\subset \cdots \subset \mathcal A_n \subset \dots$$
\item The height of $\mathcal A$,  is a limit ordinal, provided that it is defined.
\end{enumerate}
\end{lemma}
See Proposition 4.7 of \cite{pouzetsobranisa}. Notice that the converse  of (2) does not hold in general, a fact which causes some complications.

The relationship  between almost inexhaustibility and
inexhaustibility is based upon properties of reductions.

Let   $R:=(E, (\rho_{i})_{i
\in I})$ be  a relational structure and $F$ be a subset of $E$. A {\it reduction} of $R$ over $F$ is a relational  structure
$M:=(E\setminus F, (\tau_j)_{j\in J})$ such that the local automorphisms of $M$ are precisely the
restrictions of local isomorphisms of $R$ fixing $F$  pointwise.

\begin{lemma}\label {pousob}
Let $R$ be a relational structure, $K(R)$ be its kernel, $r:= \vert K(R)\vert $,    and  $M$
be a reduction of $R$ over $K(R)$.  If  $K(R)$ is finite then:

\begin{enumerate}
\item $K(M)$ is empty.
\item $H(\mathcal{A}(R))=H(\mathcal{A}(M))+p$, for some integer $p$, if and only if one of these
heights exists.
\item $\varphi_R(n)\leq 2^r\varphi_S(n)$ and $\varphi_S(n)\leq a\varphi_R(n+k)$ for some $a\in \R_+^{*}$, $k\in \N$, and all $n\in \N$.

\end{enumerate}
\end{lemma}

Item 1 is  special case of Theorem 20 \cite{gibson}. Item 2 is a special case of Theorem 4.6 of \cite{pouzetsobranisa}.
Item 3 can be proved in the same way as Theorem 21 of \cite{gibson}.

The growth of a map $\varphi: \N \rightarrow \N$ is {\it invariant under translation} if
$\varphi$ and its  translate $n \rightarrow \varphi (n+1)$  have the same
growth.  In this case,  $\varphi$ is bounded by  an exponential
function. Clearly, polynomial functions are invariant under translation. Under the hypotheses of Lemma \ref{pousob} we get from Item (3):

\begin {corollary}\label{reductiongrowth}
$\varphi_R$ is bounded by a polynomial of degree $k$, resp. has polynomial growth of degree $k$,  iff $\varphi_M$ is bounded by a polynomial of degree $k$,  resp. has polynomial growth of degree $k$.
\end{corollary}

This tells us  that in order to prove  that the growth of the profile of a  relational structure with finite kernel
is polynomial it suffices to prove that the growth of a reduction over its kernel is polynomial.

\subsection{Product of a Finite Relational Structure  by  a Chain}\label{subsectionproduct}

Let $S: = (V, (\rho_{i})_{i \in I})$ be a relational structure of signature $\mu$ and  $L: =
(D,
\leq)$ be a chain. A relational structure $R$ of signature $\mu$
is a   {\it product } of $S$ and $L$, denoted by $S
\bigotimes L$, if it satisfies the following conditions.

\begin{enumerate}
\item[{(a)}] The domain is $V\times D$.
\item[(b)]  For every  $y\in D $, the map  $x \rightarrow (x, y)$ is an isomorphism
from  $S$ into $R$.
\item[(c)] For each local isomorphism  $f$ of $L$, the map $(1_{V},
f)$ which to
 $(x,y)$ associates $(x, f(y))$ is a local isomorphism  of $R$.
\end{enumerate}

If $\vert V\vert=1$, this reduces to say that every local isomorphism of  $L$ is a local isomorphism of $R$. As we have said, a relational  structure is  {\it chainable } if there is such a chain defined on
its domain.

In this context, a relational structure is   {\it almost multichainable}, resp. {\it almost chainable}, if its kernel is finite and a reduction over its kernel is a multiple of a finite, resp. a one-element, relational structure by a chain.

A relational structure  $M$ {\it freely interprets } a relational structure  $R$ on the same set
if each local isomorphism of $M$ is a local isomorphism of  $R$.
Clearly, if a relational structure $M$ is almost multichainable, every relational structure freely interpretable by $M$ has the same property.  In fact, an almost multichainable relational structure is freely interpretable by an almost multichainable relational structure of a special form:

Let  $R:= (E,
(\rho_i)_{i \in I})$ be  an   almost multichainable relational structure;  let  $F$  and $V$ be two finite sets, $L:=(D, \leq)$ be a chain such that   $E= F\cup V\times D$  and every local isomorphism of $f$ extended by the identity on $F$ and $D$ induces a local isomorphism of $R$.  Let $r:= \vert F\vert$, $m:=\vert V\vert $, $a_1, \dots, a_r$, resp. $v_1, \dots, v_m$,  be an enumeration of the members of $F$, resp. of $V$. Set $M:= (E, A_1, \dots, A_r, V_{1},\cdots V_{m}, Eq, W)$ be the relational  structure in which $A_1, \dots, A_r$ and $V_1, \dots, V_m$ are unary relations, $A_i:= \{a_i\}$ for $i=1, \dots, r$, $V_i:= \{i\}\times D$,  $Eq$ is an equivalence relation
defined by $(u, v) \in Eq$ if either $u=v$ or $u=(x,y)$ and  $v=(x', y)$, $W$ is an order defined by $(u, v)\in W$  if either $u=v$ or $u=(x, y)$, $v=(x', y')$ and $y<y'$ in $L$.

Clearly, $M$ interprets $R$. This relational structure has a finite signature. The fact that an almost multichainable  relational structures is freely-interpretable by a relational structure with finite signature allows to use the compactness theorem of first order logic. This yields:

\begin{proposition} \label{compactness}
If a relational structure $R$ is almost multichainable then every relational structure $R'$ such that $\mathcal A(R')\subseteq \mathcal A(R)$ is also almost multichainable.
\end{proposition}

\noindent{\it Proof.\enskip}
Let $M$ as above which interprets freely $R$. If some relational
structure $R'$ verifies $\mathcal A(R')\subseteq \mathcal A(R)$ then,
since the signature of $M$ is finite, the compactness theorem yields
some relational structure $M'$ such that $\mathcal A(M') \subseteq
\mathcal A(M)$ and $M'$ interprets freely $R'$. The multirelation $M$ is
almost multichainable. It is easy to see that $M'$ is almost
multichainable too. Hence $R'$ is almost multichainable too.\hfill $\blacksquare$

\begin{theorem}\label{productwqo}
Let $\mathcal A$ be the age of an almost multichainable relational structure. Then:
\begin{enumerate}
\item $\mathcal A$ is very well-quasi-ordered;
\item  $H(\mathcal A) < \omega^{\omega}$.
\end{enumerate}
\end{theorem}

\noindent{\it Proof.\enskip}Let $R$ such that $\age(R)=\mathcal A$. Let $R'$ be a reduction  over $K(R)$. Then  $R'$ is the product  $S \bigotimes L$ of  a finite relational structure $S$ by a
chain $L$.
Let $V$ be the domain of $S$, $m : = \vert V \vert$, $D$ be the domain of $L$ and $
M$ be the multirelation interpreting $R'$.

(1) In order to prove  that  $\mathcal A$ is v.w.q.o. it suffices to prove that $\mathcal A(M)$ is v.w.q.o.  For that, we consider the alphabet
$A$ whose letters are the non-empty subsets of $V$
(that is $A  = \mathfrak{ P}(V) \setminus
\{\emptyset\}$). We order $A$ by inclusion; we extend this order to an order,  compatible with the concatenation of words, on  $A
^{\ast}$, the set  of words over $A$. With this ordering $A^{\ast}$ is isomorphic to  $\mathcal A(M)$. Thus from Theorem \ref{higman2},  $\mathcal A(M)$ is w.q.o. and in fact v.w.q.o.

(2)  Since $M$ interprets $R'$  we have $o(\age (R'))\leq o(\age (M))$. We also have $o(\mathcal A(M)) \leq o(A^{\ast})$ and,  since $\vert  A \vert = 2^{m}-1$, we have  $o(A^{\ast}) =
\omega^{\omega^{2^{m}-2}}$ from Lemma \ref {lengthalphabet}.  Thus, $o(R')\leq  \omega^{\omega^{2^{m}-2}}$. Theorem  \ref{longueur} yields $H(\mathcal A(R'))\leq
\omega^{2^{m}-1}$. From this inequality, Lemma \ref{pousob} yields $H(\mathcal A)<\omega^{2^{m}-1}+\omega$.\hfill $\blacksquare$\vskip10pt

Combining (1) of Theorem \ref{productwqo} and (3) of Theorem \ref{pouzet 73}, we get:
\begin{theorem} 
An almost multichainable relational structure having a finite signature has only finitely many bounds.
\end{theorem}
For the  special case of chainable relations  this conclusion was obtained by  C.Frasnay  in 1965 \cite{frasnay} by means of his theory of indicative groups. A truly reamarkable consequence also due to Frasnay is  that for every integer $m$ there is an integer $f(m)$ such that a monomorphic relational structure with signature bounded by $m$ and size at least $f(m)$ is chainable (see Chapter  13 of \cite{fraissetr}).

\begin{lemma}[\cite{pouzetrm}]\label{inep}
Let  $R'$ be  a relational structure and  $L$ be a chain with at least two elements. Then
the kernel of $R'$ is empty if and only if  there is a product $R'':
= R'
\bigotimes L$ such that $\mathcal A(R'') = \mathcal A(R')$.
\end{lemma}

\begin{corollary}\label{multiple1}
 The  age  $\mathcal A$ of a denumerable  relational structure with empty kernel is the union of an  infinite sequence $$\mathcal A_0 \subseteq \mathcal A_1\subseteq \cdots \subseteq \mathcal A_n\subseteq\cdots$$
where each $\mathcal A_n$ is the age of  a product of some $S\in
\mathcal A_n$ by an infinite chain.
\end{corollary}

\noindent{\it Proof.\enskip}Let $R'$ be a denumerable relational
structure such that $\age(R')=\mathcal A$ and $L$ be an infinite chain.
According to
Lemma \ref{inep}, there is a product  $R'' : = R' \bigotimes L$  such that $\age(R'')=\mathcal A$. Let  $F_{0}
\subseteq \ldots \subseteq F_{n} \subseteq \ldots$ be an non-decreasing sequence of finite subsets of the domain $E'$ of $R'$ whose union is $E'$. Let $R''_n:= R''_{\restriction F_n\times L}$ and $\mathcal A_n:= \age( R''_n)$ for $n\in \N$. The sequence
$\mathcal A_n, \ldots, \mathcal A_n, \ldots$ is non-decreasing and $\mathcal A= \bigcup_{n\in N} \mathcal A_n$.\hfill $\blacksquare$

\begin{proposition}\label{multiple2}
If the age  $\mathcal A$  of a denumerable  relational structure is inexhaustible then
\begin{enumerate}
\item[--]Either $\mathcal A$ is the age of a multichainable relational structure and $H(\mathcal A)<\omega^2$.

\item[--]Or for every integer $k$, $\mathcal A$ contains the age  $\mathcal A_k$ of  an almost multichainable relational structure such that $H(\mathcal A_k)=\omega.(k+1)$.
 \end{enumerate}
\end{proposition}

\noindent{\it Proof.\enskip}

\noindent{\bf Claim.} Either $\mathcal A$ is the age of a multichainable
relational structure or for every integer $k$ it contains the age
$\mathcal A_k$ of a multichainable relational structure such that
$H(\mathcal A_k)\geq \omega.(k+1)$.

\noindent{\it Proof of the claim.} Apply Corollary \ref{multiple1}. If the sequence $$\mathcal A_0 \subseteq \mathcal A_1\subseteq \cdots \subseteq \mathcal A_n\subseteq\cdots$$  is eventually constant, then $\mathcal A=  \mathcal A_{n_{0}}$ for some non-negative integer $n_0$. Let $S_{n_0}:= R'_{\restriction F_{n_0}}$ and $R:= R''_{n_0}=R''_{\restriction F_{n_{0}}\times L}$. Then $R$  is a product of $S_{n_0}$ by $L$.
If this sequence is not eventually constant, it contains an infinite subsequence which is strictly increasing, say:
$$\mathcal A_{n_{0}} \subset \ldots \subset \mathcal A_{n_{i}} \subset \ldots $$ Since  $R''_{n_{i}}$ is a product of  product  of $S_{n_i}$ by $L$,  its kernel is empty, that is $\mathcal A_{n_i}$ is
inexhaustible. Hence, from (1) of Lemma \ref{limit} there is a strictly increasing sequence of ages between  $\mathcal A_{n_{k-1}}$ and  $\mathcal A_{n_{k}}$ . Since $R''_{n_{i}}$ is multichainable, $\mathcal A_{n_{k}}$  has an height.   
It follows that $H(\mathcal A_{n_{k}})\geq \omega.(k+1)$.\hfill $\blacksquare$\vskip10pt

Now, suppose that $\mathcal A$ is the age of a multichainable relational structure. Then   $H(\mathcal A)$ is defined. If $H(\mathcal A)<\omega^2$, the conclusion of the proposition holds. If $H(\mathcal A)\geq \omega^2$, then for every $k$, $\mathcal A$ contains an age $\mathcal A_k$ of height $\omega.(k+1)$. According to Proposition \ref{compactness}, $\mathcal A_k$ is the age of an almost multichainable relational structure.  Suppose that $\mathcal A$ is not the age of a multichainable relational structure.  Let $k$ be an integer. According to our claim, $\mathcal A$ contains  the age  $\mathcal A_k$ of  an almost multichainable relational structure such that
 $H(\mathcal A_{k})\geq \omega.(k+1)$. Let $\mathcal A'\subseteq \mathcal A_k$ be an age of height $\omega.(k+1)$. According to Proposition \ref{compactness}, $\mathcal A$ is the age of an almost multichainable relational structure. This proves the proposition. \hfill $\blacksquare$\vskip10pt

 For relational structures with finite signature we have a similar result with a more intricate proof.
 \begin{theorem}\label{omegacarre}
Let  $\mathcal A$ be the age of infinite relational structure $R$ with finite signature.
\begin{enumerate}
\item  If $H(\mathcal A) < \omega^{2}$ then $R$ is an almost multiple of a finite relation by a chain;
in particular, $R$ has  a finite kernel.
\item  If $\mathcal A$ has no height then  $\mathcal A$ contains  some age of  height
$\omega^{2}$ (and,  in fact, ages of height $\omega^2+p$ for every integer $p$).
\end{enumerate}
\end{theorem}
The proof is given below for relational structures made of binary or unary relations;
beyond, the proof is quite involved.

\noindent{\it Proof.\enskip}
 (1) We argue by induction on the height. Let $\alpha$, $\omega \leq \alpha<\omega^2$. Suppose
that  Property (1) holds for every age $\mathcal A$ such that $\omega \leq H(\mathcal A)<\alpha$.  Let
$\mathcal A$ be an age with
$H(\mathcal A)=\alpha$. Let $R$ such that
$\mathcal A(R) = \mathcal A$. First , $K(R)$, the kernel of  $R$, is finite. Indeed, let $x \in K(R)$. We have
$\mathcal A(R_{-x})
\subsetneqq \mathcal A$. Hence, from the induction hypothesis,
$\mathcal A(R_{-x})$ is the age of an almost multiple of a finite relation by a chain. According to Theorem \ref{productwqo}, $\mathcal A(R_{-x})$ is
very well-quasi-ordered.  Since $R$ is made of binary and unary relations, it follows that
$\mathcal A(R)$ is very-well-quasi-ordered too. From  (2) of  Theorem  \ref {pouzet 73}, $K(R)$ is finite.  Let  $R' $ be a reduction of
$R$ over $K(R)$. Let $n,p<\omega $ such that
$\alpha=\omega n+p$. According to (2) of Lemma  \ref{pousob}, we have
$H(\mathcal A(R')) =
\omega n$ and  $K(R') = \emptyset$. Apply Proposition \ref{multiple2}.\\
(2)If  $\mathcal A$ has no height then $\mathcal A$ is not w.q.o., hence it contains some age $\mathcal A'$ which is w.q.o.
and has infinitely many bounds in $\mathcal A$ (Lemma \ref{nash-williams}). Being
w.q.o.,
$\mathcal A'$  has an height. If
$H(\mathcal A')<\omega^{2}$, then   from $(1)$ $\mathcal A'$ is the age of an almost multiple of a finite relation
by a chain, hence is v.w.q.o. But according to (3) of Theorem \ref {pouzet 73}, it cannot have
infinitely many bounds, a contradiction. Hence $H(\mathcal A') \geq
\omega^{2}$.\hfill $\blacksquare$


 \subsection{Profile of Almost Multichainable Relational Structures}

Let  $\mathcal A$  be the age of a product  $R:=S
\bigotimes L$ with $S:= (V, (\rho_{i})_{i \in I}) \in \mathcal A$. Let ${\mathcal J}$ be a non-empty initial segment of
$\mathfrak{ P}(V)$. Let ${\mathcal F}_{{\mathcal J}} : = \{F \in \lbrack V \times E
\rbrack^{<
\omega} : F^{-1}(y) \in {\mathcal J}$ for every $ y \in E\}$. Let $\mathcal A_{{\mathcal J}} : = \{\tau(S_{\restriction F}) : F
\in {\mathcal F}_{{\mathcal J}}\}$.

\begin{lemma}
\begin{enumerate} 
\item[{(i)}]   The age ${\mathcal A}_{{\mathcal J}}$ is inexhaustible.
\item[(ii)]   ${\mathcal A}_{{\mathcal J}} \subseteq \mathcal A_{{\mathcal J}\prime}$ if ${\mathcal J} 
\subseteq {\mathcal J}^{\prime}$;
\item[(iii)]  ${\mathcal A}_{\mathfrak {P}(V)} = {\mathcal A}$.
\end{enumerate}
\end{lemma}

Let  ${\mathcal J}$ be a non-empty initial segment of
$\mathfrak{ P}(V)$ which is  minimal w.r.t. inclusion such that $\mathcal A_{{\mathcal J}} = \mathcal A$.
Let $V'$ be  a maximal element of
 ${\mathcal J}$. Let $m' := \vert V' \vert$
and
${\mathcal J'} : = {\mathcal J} \setminus \{V'\}$.  Note that if  ${\mathcal J'}
=
\emptyset$ then  $\mathcal A$ is the age of a chainable relation, in which case its profile is $1$.

Given an integer  $p$, let  $\mathcal A_{p} : = \{\tau(S_{\restriction  F}) : F \in {\mathcal
F}_{{\mathcal J}}$ and $\vert \{i : F^{-1}(i) = V'\} \vert \leq p\}$.

\begin{lemma} \label {produitcroissant}
\begin{enumerate}
\item $\mathcal A_{p}$ is an age for
every $p<\omega$.
\item $\mathcal A = \bigcup \{\mathcal A_{p} : p < \omega\}$.
\item  $\mathcal A_{p} \neq \mathcal A$  for every $p<\omega$.
\item $\mathcal A_{ p} \subseteq \mathcal A_{p+1}$.
\item  If $\mathcal A_{p} \subsetneqq \mathcal A_{p+1}$ then $\mathcal A_{p+1} \subsetneqq \mathcal A_{p+2}$.
In other words, if $p_{0}$ is minimum such that $\mathcal A_{p_{0}} \neq\mathcal  A_{0}$ then $\mathcal A_{p}
\subsetneqq \mathcal A_{p+1}$ for every $p \geq p_{0}$.
\item   For each $p \geq p_{0}$, the kernel of any relation with age
$\mathcal A_{p}$ has
$m'.p$ elements.
\end{enumerate}
\end{lemma}

\noindent{\it Proof.\enskip}Sketch.  We may suppose that  $L$  is the chain $\Q$ of  rational numbers. In the product $V\times \Q$, we select
a subset
$X$ made of ``slices''
$F \times \{r\}$ with $F \in {\mathcal J}$ and $r\in \Q$, in such a way
that between two rational all possible slices appear. Obviously
$\mathcal A(R_{\restriction X})= \mathcal A$. Let $X_p$ obtained by
deleting from $X$ all slices of the form $V'\times \{r\}$ except $p$
such slices. Then $\mathcal A(R_{\restriction X_p})= \mathcal A_p$. To
obtain that, for $p$ larger than some non-negative integer $p_0$,
$K(R_{\restriction X_p})$ is made of these $p$ slices, apply Lemma
\ref{minusab}. \hfill $\blacksquare$

\begin{lemma} \label{binomialbound}Let  $\mathcal A$  be the age of a product  $R:=S
\bigotimes L$ with $S:= (V, (\rho_{i})_{i \in I}) \in \mathcal A$. If
$H(\mathcal A)=
\omega(k+1)$ then
 $\varphi_{R}(n) \leq {n+k \choose
k}$ for every  integer $n$.
\end{lemma}

\noindent{\it Proof.\enskip}Induction on   $k$. We apply Lemma \ref{produitcroissant}. For  $p<\omega$, we denote by $R_{p}$ a relational structure  with age
$\mathcal A_{p}$ ,   by  $R'_{p}$  a reduction  over its kernel and by $\mathcal A'_p$ its age.
We have :
$$\varphi_{R}(n) \leq \varphi_{R'_{0}}(n) + \varphi_{R'_{p_{0}}}(n - p_{0}m') +
\varphi_{R'_{p_{0}+1}}( n - (p_{0}+1)m') + \ldots$$
Set $$\varphi_{k+1} (n) : =  {k-1+n \choose n} +
{k-1+n-1 \choose n-1} + \ldots + {k-1+0 \choose 0}$$ that is  $$\varphi_{k+1} (n)={k-1+n \choose
n}+\varphi_{k} (n)$$ With Pascal identity this yields:
$$\varphi_{k+1} (n)={k+n \choose n} $$
According to (1)  of Lemma
\ref{pousob}  we have
$H(\mathcal A_{p}) = H(\mathcal A'_{p}) + n_{p}$
for some  $n_{p}<\omega$. Hence $H(\mathcal A'_{p}) \leq \omega k$. Via the inductive hypothesis,
$\varphi_{R'_{p}} (n)
\leq \varphi_{k} (n)$ for every  $n$ and every  $p$.
This yields:
$$\varphi_{R'_{0}}(n) \leq \varphi_{k}(n)$$
$$\varphi_{R'_{p_{0}}} (n-p_{0} m') \leq \varphi_{k}(n-p_{0} m')$$
$$\varphi_{R'_{p_{0}+ \ell}} (n- (p_{0} + \ell) m') \leq \varphi_{k}(n -
(p_{0} +
\ell) m')$$

$$\varphi_{R}(n)\leq  \sum^{n}_{\ell = 0} \varphi_{R'_{p_{0}+\ell}}(n-(p_{0} + \ell)
m') \leq \sum^{n}_{\ell = 0} \varphi_{k}(\ell)=\varphi_{k+1}(n) $$
\hfill $\blacksquare$


We get the following result mentioned in the introduction.

\begin{theorem}\label{theoremproved} 
If the age $\mathcal A$ of a denumerable relational structure $R$ is inexhaustible and its  height $H(\mathcal A)$ is at most
$\omega(k+1)$ then $\varphi_{R}(n) \leq {n+k \choose
k}$ for every  integer $n$.
\end{theorem}

\noindent{\it Proof.\enskip}We have $H(\mathcal A)=\omega k'+p$ with $k'\leq k+1$ and $p<\omega$. Since $\mathcal A$ is inexhaustible, (1) of Lemma \ref{pousob} tells us that $p=0$.  According to Proposition \ref{multiple2}
there is some $R$ having age  $\mathcal A$ which is a product of some $S \in
\mathcal A$ by a chain. According to Lemma \ref{binomialbound}, $\varphi_{R}(n) \leq {n+k'-1 \choose
k'-1}$ for every  integer $n$. Since  ${n+k'-1 \choose
k'-1}\leq  {n+k \choose
k}$, the desired conclusion follows.\hfill $\blacksquare$\vskip10pt

\begin{lemma}\label{counterpart}
Let  $\mathcal A$  be the age of a product  $R:=S
\bigotimes L$ with $S:= (V, (\rho_{i})_{i \in I}) \in \mathcal A$. If
$H(\mathcal A)=
\omega(k+1)$ then
$\varphi_{R}$ grows at least as a polynomial of degree $k$.\end{lemma}
For the proof of this counterpart of Lemma \ref {binomialbound}  fix an increasing sequence $(\mathcal A_i)_{i\leq k}$ of  sub-ages of $\mathcal A$ such that  $H(\mathcal A_i)=\omega.(i+1)$. Then show that  increasing sequences  $(T_i)_{i\leq k+1}$, $T_i\in \mathcal A_i\setminus \mathcal A_{i-1}$, $\vert T_{k+1}\vert=n$,  yields $O(n^k)$ distinct types. With this we obtain:

 \begin{theorem}
\label{profilepouzet3}
Let $R:= (E, (\rho_i)_{i\in I}))$ be a relational structure. If the signature of $R$ is bounded or $R$ has a finite kernel then either $R$ is almost multichainable and $H(\mathcal A(R))=\omega(k+1)+p$ -in which case   $\varphi_{R}$ grows as a polynomial of degree  $k$- or $\varphi_R$ grows faster than every polynomial.
\end{theorem}

\noindent{\it Proof.\enskip}
Suppose that the kernel $K(R)$ of $R$ is finite. Let $R'$ be a reduction
over $K(R)$. If $\mathcal A(R')$ has an height, then according to (2) of
Lemma \ref{pousob}, $\mathcal A(R)$ too and $H(\mathcal A(R))=H(\mathcal
A(R'))+p$ for some non-negative integer $p$. If $H(\mathcal A(R'))=
\omega (k+1)<\omega^2$, then $R'$ is multichainable (Proposition
\ref{multiple2}), hence $R$ is almost multichainable, and $\varphi_R'$
grows as a polynomial of degree $k$ (Lemma \ref{binomialbound} and Lemma
\ref {counterpart}). According to Corollary \ref{reductiongrowth},
$\varphi_R$ grows as a polynomial of degree $k$ too. If $\mathcal A(R')$
has no height or contains an age of height $\omega^2$, then for every
non-negative integer $k$ it contains the age of an almost multichainable
relational structure of height $\omega (k+1)$ (Proposition
\ref{multiple2}). Hence, according to Lemma \ref{counterpart},
$\varphi_{R'}$ grows faster than every polynomial and, via Corollary
\ref{reductiongrowth},  
$\varphi_R$ too. If $K(R)$ is infinite, then, with our hypothesis, the
signature is bounded. Suppose that $I$ is finite. It follows from
Theorem \ref{omegacarre} that $\mathcal A(R)$ contains ages of height
$\omega (k+1)$ for every non negative integer $k$, hence $\varphi_{R}$
grows faster than every polynomial. If $I$ is infinite, then since the
signature is bounded by some integer $m$, there is some finite subset
$I'$ such that the reduct $R^{\restriction I'}:= (E, (\rho_i)_{i\in
I'})$ has the same profile as $\varphi_R$ (Fact \ref{finiteprofile}),
thus this case reduces to the previous one. \hfill $\blacksquare$\vskip10pt

 Theorem  \ref{profilpouzet1},  Theorem \ref {profilepouzet2} and Theorem \ref{heightprofile} are immediate consequences of this result.
 In the case of a relational structure with empty kernel,
Theorem \ref{heightprofile}  has a more precise form. With  the help of Lemma \ref {binomialbound} we get
 the  following   result  announced in \cite{pouzetmontreal}:

\begin{theorem} 
If the  age of a relational structure $R$ is inexhaustible and $\varphi_R$ has polynomial growth of degree $k$ then $$\varphi_{R}(n) \leq {n+k \choose
k}$$ for every  integer $n$. 
\end{theorem}

We conclude this tour with the following:
\begin{problem}
Is the profile of a relational structure $R$ bounded by some exponential
whenever the age of $R$ is very-well-quasi-ordered under embeddability?
\end{problem}

\label{lastpage-01}

\end{document}